\documentclass[letterpaper, 9pt, conference, UKenglish, version=full]{ieeeconf}  % Comment this line out
\IEEEoverridecommandlockouts{}

\usepackage{ieeeconf-customizations}
\title{\LARGE \bf
A Poisson Jump-driven SDE Approach to Distributed Gradient Descent with Sparse Communication
}

\author{Marc Weber, John Paul Strachan and Christian Ebenbauer% <-this % stops a space
\thanks{M. Weber and C. Ebenbauer are with the Chair of Intelligent Control Systems,
        RWTH Aachen University, Aachen, Germany
        {\tt\small weber/ebenbauer@ic.rwth-aachen.de},
        J. P. Strachan is with the Peter-Gruenberg-Institute (PGI-14), Forschungszentrum Juelich GmbH,
Juelich, Germany {\tt\small j.strachan@fz-juelich.de}.
        }%
}%

\begin{document}

\maketitle
\thispagestyle{empty}
\pagestyle{empty}

%%%%%%%%%%%%%%%%%%%%%%%%%%%%%%%%%%%%%%%%%%%%%%%%%%%%%%%%%%%%%%%%%%%%%%%%%%%%%%%%
\begin{abstract}
  To bridge the gap between idealised communication models and the stochastic reality of networked systems, we introduce a framework for embedding asynchronous communication directly into algorithm dynamics using stochastic differential equations~(SDE) driven by Poisson Jumps.
  We apply this communication-aware design to the continuous-time gradient flow, yielding a distributed algorithm where updates occur via sparse Poisson events.
  Our analysis establishes communication rate bounds for asymptotic stability and, crucially, a higher, yet sparse, rate that provably any desired exponential convergence performance slower than the nominal, centralized flow.
  These theoretical results, shown for unconstrained quadratic optimisation, are validated by a numerical simulation.
\end{abstract}

\section{Introduction}\label{sec:introduction}
 Distributed optimization has emerged as a critical necessity in numerous modern applications, spanning large-scale machine learning paradigms like Federated Learning~\cite{mcmahan2017communication}, resource allocation in smart grids~\cite{molzahn2017survey}, and cooperative control in multi-agent systems and robotics~\cite{yang2019survey}.
The fundamental goal is to minimize a global objective function based on data or cost functions dispersed across multiple computational units or agents. Inherently, the resolution of this problem requires collaboration, necessitating communication among the distributed participants to share local information, such as gradients, parameter updates, or consensus states.

A major discrepancy often arises, however, between the theoretical models underlying many distributed optimization algorithms and the practical realities of their deployment.
While the literature on distributed optimization under communication constraints is vast and diverse, much of the influential literature, particularly on foundational algorithms like Distributed Gradient Descent (DGD) and ADMM, relies on assumptions of idealized communication links~\cite{nedic2009distributed,boyd2011distributed}.
These assumptions frequently include perfect synchronization, fixed-clock communication rates, negligible latency, and directed or even continuous (or infinitely fast) data transfer~\cite{mcmahan2017communication,nedic2009distributed,boyd2011distributed}.
Although recent advancements, notably in event-triggered optimization, address stochasticity and asynchronicity, these approaches often design the update rule first and then impose communication scarcity based on the state of the optimization (e.g., triggering a message when the local gradient significantly changes)~\cite{heemels2012introduction,tabuada2007event,notarstefano2019distributed,carnevale2023triggered}.
To our knowledge, strategies primarily focus on adapting algorithms to scarcity, rather than designing flows where channel limitations are the native drivers of the dynamics, especially in continuous-time settings.
These requirements remain costly in terms of energy consumption and are prone to performance degradation under real-world network congestion and do not consider or provide a-priori limitations on the communication channels.

This need for communication-aware design is amplified when considering modern computing architectures.
Low-level distributed systems, such as processor cores connected by high-speed, asynchronous interconnects (like AMD’s Infinity Fabric~\cite{AMD2022epyc}), and specialized systems like neuromorphic computing architectures, inherently feature a sparse, stochastic, and event-driven communication model~\cite{mead2002neuromorphic,davies2018loihi,leroux2025analog}.
In these systems, information exchange is a costly spike or packet rather than a continuous stream, and even state-of-the-art interconnects operate on time-scales much slower than efficient processing tiles~\cite{fpia2024}. Therefore, for an algorithm to be truly efficient and applicable across diverse distributed domains—from hardware fabrics to wide-area networks—it must be communication-aware from its inception, where the channel constraints dictate the structure of the mathematical flow.

To address this challenge, we propose a framework employed in control systems and networked stability analysis~\cite{Brockett2009,zhang2019networked}, focusing on continuous-time dynamics driven by random jumps.
Specifically, the random, sparse, and intermittent nature of physical communication channels and event-driven computing architectures—such as packet drops, irregular clock skews, or hardware-enforced intermittent updates—is naturally captured within this paradigm of SDEs driven by Poisson Jumps~\cite{brockett1977modelling,Brockett2009}, a framework related to piecewise deterministic Markov processes and hybrid systems modelling~\cite{Davis93,summers2010}.
As pioneered, for example, in the work of Brockett~\cite{brockett1977modelling} and further developed in areas like channel selection and networked control, this framework treats the communication channel not as a fixed parameter, but as a dynamic, stochastic process integrated directly into the continuous flow that dictates the timing and magnitude of information flow~\cite{farokhi2014stochastic}.
This powerful conceptual tool allows us to define mathematically precise conditions for stability and performance despite high channel uncertainty, forming the theoretical core of our algorithm design.

% Contribution
In this paper, we leverage Poisson-driven SDE modelling to rethink the design of distributed optimization algorithms, moving beyond simply analyzing the effect of noise after the algorithm is designed. Our objective is to design algorithms where the communication strategy is an integrated, efficient component of the flow.
Our primary contributions are fourfold:

Framework for Asynchronous, Sparse Design: We introduce the Poisson Jump-driven SDE modelling framework in the area of distributed optimization that uses continuous-time dynamics for the individual computation units and incorporates inter-unit communication only through stochastic Poisson-Jump events.
This framework allows for the formal analysis and design of distributed algorithms under native asynchronous, sparse, and even unstable communication dynamics.

Distributed Poisson-Jump Gradient Flow:
We illustrate this concept by taking the classic gradient descent flow~\cite{nedic2009distributed} and formulating its distributed counterpart, where local gradient information is aggregated between agents via sparse, asynchronous Poisson-Jump driven communication.
The resulting system is modeled as a set of coupled Poisson-driven SDEs.

Communication Rate for Stability: We establish a sufficient communication rate (related to the intensity of the Poisson process) under which the distributed optimization algorithm retains asymptotic stability (in the means square sense) toward the optimal solution.
This result provides a rigorous, quantifiable bound on the intermittent, random connectivity required for stability.

Communication Efficiency and desired Rates:
Crucially, we identify higher, yet still sparse range of communication rates, under which the distributed algorithm achieves any exponential convergence rate up to the one of the nominal, fully connected Gradient Descent flow.
While this results gives a formal characterization of the efficiency of the communication scheme, our experiments demonstrate that exceeding the communication rate for stability does not improve the worst-case convergence rate, thus defining a target for maximally efficient, minimal-communication algorithm design.

% Crucially, we identify a higher, yet still sparse, range of communication rates under which the distributed algorithm recovers the exponential convergence rate of the nominal, fully connected Gradient Descent flow.
% This result formally characterizes the efficiency of the communication scheme, demonstrating that exceeding this rate does not improve the worst-case convergence rate, thus defining a target for maximally efficient, minimal-communication algorithm design.

\section{Preliminaries}\label{sec:preliminaries}
 
To formalise the sparse and asynchronous communication model central to this work, we first establish the necessary mathematical framework.
We begin with a brief definition of the graph-theoretic notation used to describe the network topology, followed by the stochastic process model that governs the system dynamics.

\subsection{Graph-Theoretic Notation}
  We model a network of communicating agents as a directed graph \(\mathcal{G} = (\mathcal{V}, \mathcal{E})\).
  The set of nodes \(\mathcal{V} = \{\nu_1, ..., \nu_n\}\) represents the \(n\) distributed computational units or agents.
  The set of edges \(\mathcal{E} \subset \mathcal{V} \times \mathcal{V}\) defines the communication links, where an ordered pair \((j, i) \in \mathcal{E}\) denotes a directed channel from node $\nu_j$ to node $\nu_i$.
  Consequently, $\nu_j$ is an in-neighbour of $\nu_i$, and the set of all such in-neighbours is $\mathcal{N}_i = \{\nu_j \in \mathcal{V} \mid (j, i) \in \mathcal{E}\}$

\subsection{Poisson Jump-driven SDEs}

  We model the state evolution of agents and channels using SDEs driven by Poisson jumps~\cite{Brockett2009,Davis93}.
  This approach naturally captures the sporadic, event-driven nature of communication.

  % \begin{definition}[Poisson Jump-Driven SDE]\label{def:poisson-jump-driven-sde}
% \mbox{}
% Let \(\vec{x}: T \times \Omega \to \mathbb{R}^n\), written as \(\vec{x}(t, \omega)\), be a stochastic process defined on the probability space \(\parens*{\Omega, \mathcal{F}, P}\), taking values in \(\parens*{\mathbb{R}^n, \mathcal{B}\parens*{\mathbb{R}^n}}\), where \(\mathcal{B}\parens*{\mathbb{R}^n}\) is the Borel~\(\sigma\)-algebra.
% The process evolves according to the stochastic differential equation

\begin{figure}[h]
    \centering
    \includegraphics{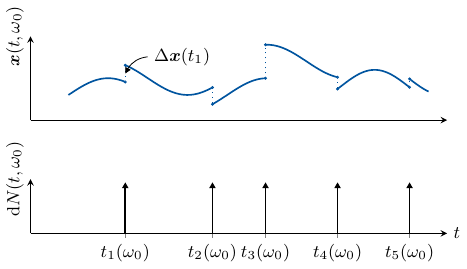}
    \caption{Sample path of the Poisson jump-driven SDE, showing the continuous evolution of \(\vec{x}(t, \omega_0)\) interrupted by discrete jumps.}\label{fig:example-trajectory-of-sde}
  \end{figure}

Consider a complete probability space~\(\parens*{\Omega, \mathcal{F}, P}\), equipped with a filtration~\(\braces{\mathcal{F}_t}_{t \geq 0}\) satisfying the usual conditions.
We model the system state as a \(\braces{\mathcal{F}_t}\)-adapted stochastic process~\(\vec{x}(t, \omega): T \times \Omega \to \mathbb{R}^n\), taking values in \(\R^n\).
Consistent with the definition of the solutions in the Itô sense, we write its dynamics in differential form as
\begin{equation}
  \label{eqn:SDE-prototype}
  \odif{\vec{x}}(t, \omega) = f\parens*{\vec{x}(t, \omega)} \odif{t} + g\parens*{\vec{x}(t, \omega)} \odif{N}(t, \omega),
\end{equation}
where \({f: \mathbb{R}^n \to \mathbb{R}^n}\) is the continuous drift vectorfield and~\({g: \mathbb{R}^n \to \mathbb{R}^n}\) is the jump map.
The term~\(N(t, \omega)\) represents a standard \(\braces{\mathcal{F}_t}\)-Poisson counting process with rate~\(\lambda\).
This model~\eqref{eqn:SDE-prototype} yields a sample path~\(\vec{x}(t, \omega_0)\) for a fixed outcome~\(\omega_0 \in \Omega\), that is cádlág, i.\,e.\ right-continous with left limits, which is denoted as
\begin{equation*}
    \vec{x}(t, \omega_0) \!= \!
    \vec{x}(t_0, \omega_0) 
    + \!\!\int_0^t \!\!\!\!f(\vec{x}(\tau, \omega_0)) \odif{\tau}
    + \!\!\int_0^t \!\!\!\! g(\vec{x}(\tau, \omega_0)) \odif{N}(\tau, \omega_0).
\end{equation*}
% The notation~\({\vec{x}(t^-, \omega) = \lim_{s \nearrow t} \vec{x}(s, \omega)}\) denotes the left-limit of the state, ensuring the system state just before a potential jump is used to calculate the jump magnitude.
The interpretation of these pathwise dynamics, as illustrated in \Cref{fig:example-trajectory-of-sde}, 
consists of two key behaviours (see \cite{Brockett2009,Davis93} for details):
\begin{description}
  \item[Continuous Flow]
        Between jump (event) times of the Poisson process \(N(t, \omega)\), the system evolves deterministically according to the ordinary differential equation~\(\dot{\vec{x}}(t) = f(\vec{x}(t))\).
  \item[Discrete Jumps]
        At each jump time~\(t_i(\omega_0)\), the state experiences an instantaneous, discrete shift~\(\Delta \vec{x}(t_i, \omega_0) :=\vec{x}(t_i, \omega_0)-\vec{x}(t_i^-, \omega_0) := g(\vec{x}(t_i^-, \omega_0))\) where~\({\vec{x}(t^-, \omega) = \lim_{s \nearrow t} \vec{x}(s, \omega)}\) denotes the left-limit of the state, representing the value immediately preceding the jump.
\end{description}

% is the differential of a Poisson process~\(N(t, \omega)\) with rate \(\lambda\), taking the value 1 at the jump times \(\braces*{t_i(\omega)}\) and 0 otherwise, with \({P\parens*{\odif{N}(t, \omega) = 1} \approx \lambda \odif{t}}\) over an infinitesimal \(\odif{t}\).
% A sample path \(\vec{x}(t, \omega)\) for a fixed \(\omega \in \Omega\) is piecewise right-continuous: it evolves continuously via \(f\) between jump times and, at each \(t_i(\omega)\), jumps by \(\Delta \vec{x}(t_i, \omega) = g(\vec{x}(t_i^-, \omega))\), where \({\vec{x}(t_i^-, \omega) = \lim_{s \nearrow t_i} \vec{x}(s, \omega)}\) is the left limit, ensuring right-continuity at \(t_i\).
% % \end{definition}

  For a function~\(V: \mathbb{R}^n \to \mathbb{R}\), its evolution under~\eqref{eqn:SDE-prototype} is governed by the \emph{Itô Differentiation Rule}
  \begin{multline}
    \odif{V}(\vec{x}(t, \omega)) =
    % \\
    % \begin{aligned}
    \pdv{V}{\vec{x}}(\vec{x}(t, \omega)) f\parens*{\vec{x}(t, \omega)} \odif{t}
    \\
    + \brackets*{V(\vec{x}(t, \omega) + g\parens*{\vec{x}(t, \omega)}) - V(\vec{x}(t, \omega))} \odif{N}(t, \omega).
    % \end{aligned}
  \end{multline}
  For the stability analysis in \Cref{sec:main-results}, we move from the pathwise Itô rule to the expected evolution of a Lyapunov function candidate~\(V(\vec{x})\).
  %This is governed by the infinitesimal generator \(\mathcal{L}\) of the process.
  For a sufficiently smooth function~\({V: \mathbb{R}^n \to \mathbb{R}}\), this is given by
  \({\odv*{\E{ V(\vec{x}(t, \omega))}}{t} =  \E{\mathcal{L}V(\vec{x}(t, \omega))}}\)  with
  \begin{equation*}
    \mathcal{L}V(x) = \pdv{V}{\vec{x}}(\vec{x})  f(\vec{x})
    + \lambda \brackets{
        V(\vec{x} + g(\vec{x})) - V(\vec{x})
      }.
  \end{equation*}

  This operator describes the expected rate of change of \(V(\vec{x})\) of a sample path of \eqref{eqn:SDE-prototype} passing through $\vec{x}$.
  This forms the basis for the Lyapunov stability analysis via~\Cref{lemma:gpues-ms-lyapunov-like-theorem}.

\subsection{Stability Concepts}
  The following definition and theorem are adapted from~\cite[Definition~4.1, Theorem~4.4]{mao2007stochastic} with an extension for practical stability taken from the concept of ISS-stability~\cite[Theorem~3.4]{sontag2008input} for constant inputs to characterize the stability properties obtained in this work.

  %  the formal definition of practical uniform exponential stability in the measn-square sense
  % TODO optional: make brackets of expectation operator bigger?
% TODO optional: use interval notation in last line?
\begin{definition}[Global Practical Uniform Exponential Stability in mean-square sense]\label{def:gpues-ms}
  A point \({\vec{x}^* \in \R^n}\) is said to be \emph{globally practically uniformly exponentially stable in the mean square sense} for a SDE  system~\eqref{eqn:SDE-prototype}, if there exist constants~\({\alpha, \beta \in \interval[open]{0}{\infty}}\) and \(\gamma \in \interval[open right]{0}{\infty}\) such that for any initial condition \(\vec{x}(t_0) = \vec{x}^0 \in \R^{n}\), the  inequality
  \begin{equation}
    \label{eqn:gpues-ms-lyapunov-condition}
    \mathbb{E}[\norm{\vec{x}(t) - \vec{x}^*}^2] \le \alpha \, \norm*{\vec{x}^0 - \vec{x}^*}^2 e^{-\beta(t-t_0)} + \gamma
  \end{equation}
  holds for all \(t \in \interval[open right]{t_0}{\infty}\).
\end{definition}
  If \(\gamma = 0\), then the attribute \emph{practically} can be dropped.
  We characterize this stability concept by means of a Lyapunov-like theorem.

  % the Lyapunov-like characterization of ms-pues
  
% TODO optional: make enumeration symbols smaller
\begin{lemma}[Global Practical Uniform Exponential Stability in mean-square sense]\label{lemma:gpues-ms-lyapunov-like-theorem}
  Consider the SDE with jumps~\eqref{eqn:SDE-prototype}.
  Assume there exists a function \({V: \mathbb{R}^n \to \mathbb{R}}\) and positive constants \({c_1, c_2, c_3, \gamma'}\) such that
  \begin{enumerate}[
      label={\tiny(\arabic*)},
      wide,
      % labelindent=0pt,
      labelsep=0.75em
    ]
    \item \(c_1 \norm{\vec{x} - \vec{x}^*}^2 \le V(\vec{x}) \le c_2 \norm{\vec{x} - \vec{x}^*}^2\) for all \(\vec{x} \in \mathbb{R}^n\) and
    \item \(\mathcal{L}V(\vec{x}) \le -c_3 \norm{\vec{x} - \vec{x}^*}^2 + \gamma'\) for all \(\vec{x} \in \mathbb{R}^n\).
  \end{enumerate}
  Then, the state \(\vec{x}^*\) is \emph{globally practically uniformly exponentially stable in the mean square sense} for~\eqref{eqn:SDE-prototype}, with parameters~\(\alpha = \frac{c_2}{c_1} \), \(\beta = \frac{c_3}{c_2}\) and \(\gamma = \frac{\gamma'}{c_3}\).
\end{lemma}

\section{Main Results}\label{sec:main-results}
 % This file contains the system model, the Lyapunov function candidate,
% and the full derivation of its infinitesimal generator.

In this section, we apply the previously established SDE framework to bridge the gap between idealized communication models often assumed in distributed optimization and the sparse, asynchronous reality of networked systems.
Our goal is to design a distributed gradient descent algorithm where communication constraints are integral to the dynamics.

\subsection{Channel Dynamics Model}\label{sec:channel-dynamics}

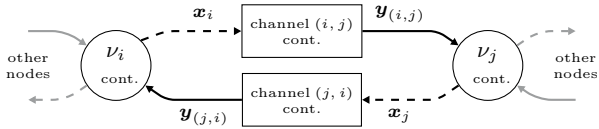
\begin{figure}[hbt]
    \centering
    \begin{tikzpicture}[
  font=\small,
  % 4. Radius als Teil des Stils
  edge radius/.initial={5pt},
  node style/.style={
      circle,
      draw,
      align=center,
      text width=4mm,
      font=\small,
    },
  block/.style={
      rectangle,
      draw,
      minimum width=1.5cm,
      minimum height=0.7cm,
      align=center,
      font=\tiny,
    },
  % 3. Kleinere, einfache Dreiecke (z.B. Classic, scale=0.7)
  edge style/.style={
  ->,
  -{Triangle[width=1mm,length=1mm]},
  thick,
  rounded corners=3pt,
  },
  edge jump style/.style={
  ->,
  -{Triangle[width=1mm,length=1mm]},
  thick,
  dashed,
  rounded corners=3pt,
  },
  % 1. Angedeutete Kanten: Dick, Black50, Dashed, kleiner Pfeil
  dummy edge/.style={
  ->,
  -{Triangle[width=1mm,length=1mm]},
  thick,
  black50,
  % dotted, % Hier wieder dashed eingefügt, da es vorher entfernt wurde
  rounded corners=\pgfkeysvalueof{/tikz/edge radius}
  },
  dummy edge jump/.style={
  ->,
  -{Triangle[width=1mm,length=1mm]},
  thick,
  black50,
  dashed,
  rounded corners=\pgfkeysvalueof{/tikz/edge radius}
  }
  ]

  % 1. Knoten definieren
  \node[node style] (nu_i) {\( \nu_i \)\\ {\tiny cont.}};
  \node[node style, right=4cm of nu_i] (nu_j) {\( \nu_j \)\\ {\tiny cont.}};

  % --- Kante nu_i -> nu_j (System S_ij) ---
  \node[block, above=0.1cm of $(nu_i)!0.5!(nu_j)$] (block_ij) {channel \( (i, j) \)\\ {\tiny cont.}};

  % BERECHNUNG DES KNICKPUNKTES FÜR S_ij (25% des Weges von nu_i)
  \coordinate (mid_x_ij) at ($(nu_i.north east)!0.25!(block_ij.west)$);
  \coordinate (mid_point_ij) at (mid_x_ij |- block_ij.center);

  % Verbindung nu_i -> S_ij
  \draw[edge jump style]
  (nu_i.north east) --
  (mid_point_ij) --
  node[midway, above, font=\scriptsize] {\( \vec{x}_{i\vphantom{j}} \)} % Signalname 1 (Input zum Block)
  (block_ij);

  % BERECHNUNG DES KNICKPUNKTES ZURÜCK NACH nu_j (75% des Weges von Block)
  \coordinate (mid_x_j) at ($(block_ij.east)!0.75!(nu_j.north west)$);
  \coordinate (mid_point_j) at (mid_x_j |- block_ij.center);

  % Verbindung S_ij -> nu_j
  \draw[edge style]
    (block_ij.east) --
    node[midway, above, font=\scriptsize] {\( \vec{y}_{(i, j)} \)} % Signalname 2 (Output vom Block)
    (mid_point_j) --
    (nu_j);

  % --- Kante nu_j -> nu_i (System S_ji) ---
  \node[block, below=0.1cm of $(nu_i)!0.5!(nu_j)$] (block_ji) {channel \( (j, i)\)\\ {\tiny cont.}};

  % BERECHNUNG DES KNICKPUNKTES FÜR S_ji (25% des Weges von nu_j)
  \coordinate (mid_x_ji) at ($(nu_j.south west)!0.25!(block_ji.east)$);
  \coordinate (mid_point_ji) at (mid_x_ji |- block_ji.center);

  % Verbindung nu_j -> S_ji
  \draw[edge jump style]
    (nu_j) --
    (mid_point_ji) --
    node[midway, below, font=\scriptsize] {$\vec{x}_{j}$} % Signalname 3
    (block_ji);

  % BERECHNUNG DES KNICKPUNKTES ZURÜCK NACH nu_i (75% des Weges von Block)
  \coordinate (mid_x_i) at ($(block_ji.west)!0.75!(nu_i.south east)$);
  \coordinate (mid_point_i) at (mid_x_i |- block_ji.center);

  % Verbindung S_ji -> nu_i
  \draw[edge style]
    (block_ji) --
    node[midway, below, font=\scriptsize] {\( \vec{y}_{(j, i)} \)} % Signalname 4
    (mid_point_i) --
    (nu_i);

  % 4. Angedeutete Kanten (dummy edges)
  % 1. Nur ein Paar links und ein Paar rechts, mit Knicken

  \pgfmathsetmacro{\dummyoffset}{0.5}

  % --- Angedeutete Kanten für nu_i (Links, außerhalb des Hauptdiagramms) ---

  % Von anderen Knoten (z.B. nu_k) zu nu_i
  \coordinate (tmp_i_in) at ($ 2*(nu_i) - (mid_point_i) $);
  \coordinate (dummy_mid_i_in) at (tmp_i_in |- mid_point_ij);
  \coordinate (dummy_start_i_in) at ($(dummy_mid_i_in) - (\dummyoffset, 0) $);
  \draw[dummy edge]
    (dummy_start_i_in) --
    (dummy_mid_i_in) --
    (nu_i);

  \coordinate (tmp_i_out) at ($ 2*(nu_i) - (mid_point_ij) $);
  \coordinate (dummy_mid_i_out) at (tmp_i_out  |- mid_point_i);
  \coordinate (dummy_end_i_out) at ($(dummy_mid_i_out) - (\dummyoffset, 0)$);
  \draw[dummy edge jump]
    (nu_i) --
    (dummy_mid_i_out) --
    (dummy_end_i_out);

  % --- Angedeutete Kanten für nu_j (Rechts, außerhalb des Hauptdiagramms) ---

  % Von anderen Knoten zu nu_j
  \coordinate (tmp_j_in) at ($ 2*(nu_j) - (mid_point_j) $);
  \coordinate (dummy_mid_j_in) at (tmp_j_in |- mid_point_ji);
  \coordinate (dummy_start_j_in) at ($(dummy_mid_j_in) + (\dummyoffset, 0) $);
  \draw[dummy edge]
    (dummy_start_j_in) --
    (dummy_mid_j_in) --
    (nu_j);

  % Von nu_j zu anderen Knoten
  \coordinate (tmp_j_out) at ($ 2*(nu_j) - (mid_point_ji) $);
  \coordinate (dummy_mid_j_out) at (tmp_j_out  |- mid_point_j);
  \coordinate (dummy_end_j_out) at ($(dummy_mid_j_out) + (\dummyoffset, 0)$);
  \draw[dummy edge jump]
    (nu_j) --
    (dummy_mid_j_out) --
    (dummy_end_j_out);

  % Definieren Sie Ihren Text (ersetzen Sie "Weitere Knoten" durch Ihren gewünschten Text)
  \node[font=\tiny, align=center] at ($(dummy_start_j_in)!0.5!(dummy_end_j_out)$) {
      other\\
      nodes
    };
  \node[font=\tiny, align=center] at ($(dummy_start_i_in)!0.5!(dummy_end_i_out)$) {
      other\\
      nodes
    };

\end{tikzpicture}
    \caption{
        Illustration of two nodes within the communication graph and their directed communication channels: solid lines represent continuous signals, while dashed lines represent discrete communication events. Channel and node dynamics are continuous-time dynamics. }%
    \label{figure:graph-illustration}
  \end{figure}

  To realise distributed algorithms within our SDE framework, we model the interaction between any two agents \( \nu_j \) and \( \nu_i \) connected by an edge \( (j,i) \in \mathcal{E} \).
  As depicted in \Cref{figure:graph-illustration}, each agent \( \nu_i \) maintains its own state \( \vec{x}_i \) governed by continuous-time dynamics.
  Information exchange between agents, however, occurs only through discrete, stochastic communication events modeled as Poisson jumps.
  The continuous-time evolution of the agent dynamics can be interpreted as a time scale separation in the sense that sampling rates of the agents are much faster than the ones over the communication channel.

  Specifically, for each directed communication link \( (j,i) \in \mathcal{E} \), agent \( \nu_i \) maintains a local channel state \( \vec{z}_{(j,i)}(t,\omega) \in \mathbb{R}^{d_{j}} \).
  This state represents the information agent \( \nu_i \) possesses about agent \( \nu_j \)'s state \( \vec{x}_j \).
  The crucial aspect is that \( \vec{z}_{(j,i)} \) is only updated sporadically when a communication event (a Poisson jump) occurs on the channel \( (j,i) \).
  Between these events, the channel state might evolve continuously (e.\,g., decay) or remain constant.
  The general dynamics for the channel state \( \vec{z}_{(j,i)} \) and its output \( \vec{y}_{(j,i)} \) (which is then used by agent \( \nu_i \)) are given by the SDE
  \begin{subequations}
    \begin{align}
      \label{eqn:channel_state}
      \odif{\vec{z}}_{(j, i)}(t, \omega)
                                  & =
      \begin{aligned}[t]
         & f(\vec{z}_{(j, i)}(t, \omega)) \odif{t}
        \\
         & + g(\vec{x}_j(t, \omega), \vec{z}_{(j, i)}(t, \omega)) \odif{N}_{(j, i)}(t, \omega),
      \end{aligned}
      \\
      \label{eqn:channel_output}
      \vec{y}_{(j, i)}(t, \omega) & = h(\vec{z}_{(j, i)}(t, \omega)).
    \end{align}
  \end{subequations}

  Here, \( \odif{N}_{(j,i)} \) is the Poisson process modelling communication events from \( \nu_j \) to \( \nu_i \) with average rate \( \lambda_{(j, i)} \).
  The map~\({f: \mathbb{R}^{d_{(j,i)}} \to \mathbb{R}^{d_{(j,i)}}}\) describes the continuous evolution of the channel state between communications, while~\({g: \mathbb{R}^{d_j} \times \mathbb{R}^{d_{(j, i)}} \to \mathbb{R}^{d_{(j, i)}}}\) defines how the state \( \vec{x}_j \) of the sender updates the channel state \( \vec{z}_{(j,i)} \) at the instant of a communication event.
  The output map~\(h: \mathbb{R}^{d_{(j, i)}} \to \mathbb{R}^{m_{(j,i)}}\) determines what information \( \vec{y}_{(j,i)} \) is actually available to the receiving agent \( \nu_i \) based on the channel state.

  This model explicitly separates the continuous internal dynamics of each agent from the discrete, event-driven communication between them.
  It provides flexibility to represent various channel behaviours, including potentially unstable dynamics between updates, i.\,e., \(f\) could represent an unstable system, allowing for modelling scenarios like accumulating errors or decaying information quality.

  We explore two representative effects of Poisson jumps on the channel output \(\vec{y}_{(j, i)}\), illustrating the modelling approach through practical scenarios contrasted in \Cref{fig:example-trajectory-comparison}

  \subsubsection{Stochastic Ideal Sample-and-Hold}
  In environmental monitoring networks, sensors transmit measurements like temperature or humidity to a central hub, which must retain the latest vector of readings exactly until a new message arrives.
  This requires a channel that captures \(\nu_j\)'s multi-dimensional state instantly and holds it without degradation, akin to an ideal sample-and-hold process in signal processing.

  \begin{example}[Ideal Sample-and-Hold Channel]
    Consider a channel syncing \(\nu_i\)'s data with \(\nu_j\)'s state at message arrivals.
    The dynamics are
    \begin{subequations}
      \begin{align}
        \label{eqn:ideal_sah_state}
        \odif{\vec{z}}_{(j, i)} & = (\vec{x}_j - \vec{z}_{(j, i)}) \odif{N}_{(j, i)},
        \\
        \label{eqn:ideal_sah_output}
        \vec{y}_{(j, i)}        & = \vec{z}_{(j, i)},
      \end{align}
    \end{subequations}
    where \(\vec{z}_{(j, i)} \in \mathbb{R}^{d_{j}}\).
    At each jump, \(\vec{z}_{(j, i)}\) resets to \(\vec{x}_j(t^-)\), delivering \(\nu_j\)'s state to \(\nu_i\) without loss.
  \end{example}

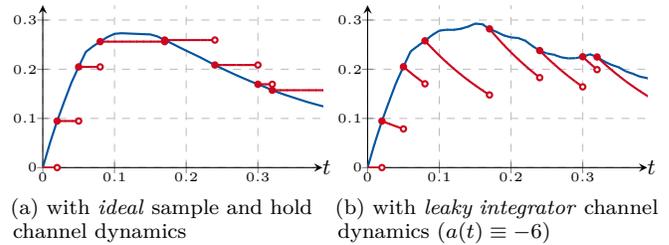
\begin{figure}[tb]
    \begin{subfigure}[t]{0.5\linewidth}%
      % --- THE OPTIMIZATION: Read the table ONCE ---
% Read the data from 'simulation_data.csv' and store it in a macro named '\data'

\pgfplotstableread[col sep=comma]{src/figures/figure.lyapunov-function-evolution/l27a0/dgd_sample_path_run_1.csv}\dataSPA
% ---------------------------------------------

\begin{tikzpicture}[
    % trim axis left,
    outer sep = 0pt,
    inner sep=0pt,
  ]

  % 1. DEFINE THE FIXED MARGIN NODE (Must be OUTSIDE the axis environment)
  % This node ensures the anchor 'margin' is known.
  \node[
      anchor=west,           % Node's (0,0) reference is its left edge
      inner sep=0pt,
      outer sep=0pt,
      minimum width=3mm,    % <--- SET YOUR FIXED MARGIN HERE
      draw=none
    ] (margin) at (0,0) {};
  %-----------------------------------------------------------------
  \begin{axis}[
      scale only axis,
      at={(margin.east)},
      anchor=south west,
      width = \linewidth-3mm-3mm,
      height = 0.5*\linewidth,
      xlabel={\(t\)},
      % ylabel={\(z(t)\)},
      % ymode = log,
      grid=major,
      xmin = 0,
      xmax = 0.39,
      ymin=0,
      ymax=0.33,
      %   legend pos=south east,
      axis lines=left,
      ylabel near ticks,
      % ylabel style={at={(ticklabel cs:1)},anchor=south,rotate=-90}, % Move y-label to the top and prevent rotation
      tick label style={font=\tiny},
      label style={font=\small},
    ]

    % 1. Draw the LINES ONLY using the custom colormap.
    \addplot [
        % mesh,
        % line legend,
        % shader=flat corner, % Use mesh to style each segment
        % thick,
        % colormap name=jumpmap_magenta_fast, % Dies ist der entscheidende Befehl!
        color = blue,
        % scatter,
        % scatter src=explicit symbolic,
        % mark size=1pt,
        % scatter/classes={
        %     0={mark=none},
        %     1={mark=*, draw=blue, fill=white},
        %     2={mark=*, draw=blue, fill=blue}
        %   },
        % opacity = {
        %     ifthenelse(\pgfplotspointmetatransformed == 1 , 0, 1)
        %   },
      ] table [
        x=time,
        y=x3_dim1,
        meta=x3_dim1_style,
        col sep=comma
      ] {\dataSPA};\label{plot:figure:samplepath-nodrift:x}

    % 1. Draw the LINES ONLY using the custom colormap.
    \addplot [
        mesh,
        line legend,
        shader=flat corner, % Use mesh to style each segment
        % thick,
        % colormap name=jumpmap_magenta_fast, % Dies ist der entscheidende Befehl!
        color = red,
        scatter,
        scatter src=explicit symbolic,
        mark size=1pt,
        scatter/classes={
            0={mark=none},
            1={mark=*, draw=red, fill=white},
            2={mark=*, draw=red, fill=red}
          },
        opacity = {
            ifthenelse(\pgfplotspointmetatransformed == 1 , 0, 1)
          },
      ] table [
        x=time,
        y=z23_dim1,
        meta=z23_dim1_style,
        col sep=comma
      ] {\dataSPA};\label{plot:figure:samplepath-nodrift:z}

    % \addplot [
    %     mesh,
    %     shader=flat corner, % Use mesh to style each segment
    %     % thick,
    %     % colormap name=jumpmap_magenta_fast, % Dies ist der entscheidende Befehl!
    %     color = maygreen,
    %     scatter,
    %     scatter src=explicit symbolic,
    %     mark size=1pt,
    %     scatter/classes={
    %         0={mark=none},
    %         1={mark=*, draw=maygreen, fill=white},
    %         2={mark=*, draw=maygreen, fill=maygreen}
    %       },
    %     opacity = {
    %         ifthenelse(\pgfplotspointmetatransformed == 1 , 0, 1)
    %       },
    %   ] table [
    %     x=time,
    %     y=V,
    %     meta=V_style,
    %     col sep=comma
    %   ] {\dataSPB};

    % \addplot [
    %     mesh,
    %     shader=flat corner, % Use mesh to style each segment
    %     % thick,
    %     colormap name=jumpmap_magenta_fast, % Dies ist der entscheidende Befehl!
    %     color = petrol,
    %     scatter,
    %     scatter src=explicit symbolic,
    %     mark size=1pt,
    %     scatter/classes={
    %         0={mark=none},
    %         1={mark=*, draw=petrol, fill=white},
    %         2={mark=*, draw=petrol, fill=petrol}
    %       },
    %     opacity = {
    %         ifthenelse(\pgfplotspointmetatransformed == 1 , 0, 1)
    %       },
    %   ] table [
    %     x=time,
    %     y=V,
    %     meta=V_style,
    %     col sep=comma
    %   ] {\dataSPC};

  \end{axis}
\end{tikzpicture}
      \caption{with \emph{ideal} sample and hold channel dynamics}
      \label{example-trajectory-comparison:a}
    \end{subfigure}%
    \begin{subfigure}[t]{0.5\linewidth}
      % --- THE OPTIMIZATION: Read the table ONCE ---
% Read the data from 'simulation_data.csv' and store it in a macro named '\data'

\pgfplotstableread[col sep=comma]{src/figures/figure.lyapunov-function-evolution/l26a-6/dgd_sample_path_run_1.csv}\dataSPA
% ---------------------------------------------

\begin{tikzpicture}[
    % trim axis left,
    outer sep = 0pt,
    inner sep=0pt,
  ]

  % 1. DEFINE THE FIXED MARGIN NODE (Must be OUTSIDE the axis environment)
  % This node ensures the anchor 'margin' is known.
  \node[
      anchor=west,           % Node's (0,0) reference is its left edge
      inner sep=0pt,
      outer sep=0pt,
      minimum width=3mm,    % <--- SET YOUR FIXED MARGIN HERE
      draw=none
    ] (margin) at (0,0) {};
  %-----------------------------------------------------------------
  \begin{axis}[
      scale only axis,
      at={(margin.east)},
      anchor=south west,
      width = \linewidth-3mm-3mm,
      height = 0.5*\linewidth,
      xlabel={\(t\)},
      % ylabel={\(z(t)\)},
      % ymode = log,
      grid=major,
      xmin = 0,
      xmax = 0.39,
      ymin=0,
      ymax=0.33,
      %   legend pos=south east,
      axis lines=left,
      ylabel near ticks,
      % ylabel style={at={(ticklabel cs:1)},anchor=south,rotate=-90}, % Move y-label to the top and prevent rotation
      tick label style={font=\tiny},
      label style={font=\small},
    ]

    % 1. Draw the LINES ONLY using the custom colormap.
    \addplot [
        mesh,
        line legend,
        shader=flat corner, % Use mesh to style each segment
        % thick,
        % colormap name=jumpmap_magenta_fast, % Dies ist der entscheidende Befehl!
        color = blue,
        scatter,
        scatter src=explicit symbolic,
        mark size=1pt,
        scatter/classes={
            0={mark=none},
            1={mark=*, draw=blue, fill=white},
            2={mark=*, draw=blue, fill=blue}
          },
        opacity = {
            ifthenelse(\pgfplotspointmetatransformed == 1 , 0, 1)
          },
      ] table [
        x=time,
        y=x3_dim1,
        meta=x3_dim1_style,
        col sep=comma
      ] {\dataSPA};\label{plot:figure:samplepath-drift:x}

    % 1. Draw the LINES ONLY using the custom colormap.
    \addplot [
        mesh,
        line legend,
        shader=flat corner, % Use mesh to style each segment
        % thick,
        % colormap name=jumpmap_magenta_fast, % Dies ist der entscheidende Befehl!
        color = red,
        scatter,
        scatter src=explicit symbolic,
        mark size=1pt,
        scatter/classes={
            0={mark=none},
            1={mark=*, draw=red, fill=white},
            2={mark=*, draw=red, fill=red}
          },
        opacity = {
            ifthenelse(\pgfplotspointmetatransformed == 1 , 0, 1)
          },
      ] table [
        x=time,
        y=z23_dim1,
        meta=z23_dim1_style,
        col sep=comma
      ] {\dataSPA};\label{plot:figure:samplepath-drift:z}

    % \addplot [
    %     mesh,
    %     shader=flat corner, % Use mesh to style each segment
    %     % thick,
    %     % colormap name=jumpmap_magenta_fast, % Dies ist der entscheidende Befehl!
    %     color = maygreen,
    %     scatter,
    %     scatter src=explicit symbolic,
    %     mark size=1pt,
    %     scatter/classes={
    %         0={mark=none},
    %         1={mark=*, draw=maygreen, fill=white},
    %         2={mark=*, draw=maygreen, fill=maygreen}
    %       },
    %     opacity = {
    %         ifthenelse(\pgfplotspointmetatransformed == 1 , 0, 1)
    %       },
    %   ] table [
    %     x=time,
    %     y=V,
    %     meta=V_style,
    %     col sep=comma
    %   ] {\dataSPB};

    % \addplot [
    %     mesh,
    %     shader=flat corner, % Use mesh to style each segment
    %     % thick,
    %     colormap name=jumpmap_magenta_fast, % Dies ist der entscheidende Befehl!
    %     color = petrol,
    %     scatter,
    %     scatter src=explicit symbolic,
    %     mark size=1pt,
    %     scatter/classes={
    %         0={mark=none},
    %         1={mark=*, draw=petrol, fill=white},
    %         2={mark=*, draw=petrol, fill=petrol}
    %       },
    %     opacity = {
    %         ifthenelse(\pgfplotspointmetatransformed == 1 , 0, 1)
    %       },
    %   ] table [
    %     x=time,
    %     y=V,
    %     meta=V_style,
    %     col sep=comma
    %   ] {\dataSPC};

  \end{axis}
\end{tikzpicture}
      \caption{with \emph{leaky integrator} channel dynamics~(\(a(t) \equiv -6\))}
      \label{fig:example-trajectory-comparison:b}
    \end{subfigure}%
    % \centering
    \caption{Sample paths with different channel dynamics, showing the continuous evolution of (\(\vec{x}_1(t, \omega)\),\ref{plot:figure:samplepath-nodrift:x}) and (\( \vec{z}_{(1, 2)}(t, \omega)\),\ref{plot:figure:samplepath-nodrift:z}) interrupted by discrete communication events \(\odif{N}_{(1,2)}(t, \omega)\).}\label{fig:example-trajectory-comparison}
  \end{figure}

  As seen in \Cref{example-trajectory-comparison:a}, the channel output \(\vec{z}_{(j,i)}\) perfectly holds the value of \(\vec{x}_j\)  sampled at the last jump time, until the next jump occurs, while \(\vec{x}_j\) proceeds to evolve in continuous-time.

  \subsubsection{Leaky Integrator Channel}
\begin{figure}[h]
    \centering
    \includegraphics{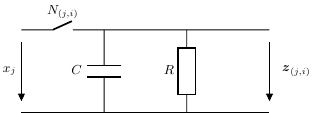}
    \caption{Switched RC circuit modelling a leaky integrator channel, where input \(x_j\) charges the capacitor at Poisson jumps, and output \(\vec{z}_{(j, i)}\) leaks through \(R\).}%
    \label{fig:rc_circuit}
  \end{figure}

  In low-cost sensor hardware, channels aim to store values but suffer from gradual loss, as in a switched RC~circuit where a capacitor leaks charge through non-ideal resistive effects, as shown in~\Cref{fig:rc_circuit}. More ideal signal holding behavior is possible, but comes at great cost in complex refresh circuits or the addition of non-volatile memories. The same effect can be observed in neurons, where activation abates over time.
  This smooths incoming messages over time, resembling a faulty sample-and-hold process.

  \begin{example}[Leaky Integrator Channel]
    Consider a channel mimicking an RC circuit. The dynamics are
    \begin{subequations}
      \begin{align}
        \odif{\vec{z}}_{(j, i)} & = a_{(j,i)} \vec{z}_{(j, i)} \odif{t} + (\vec{x}_j - \vec{z}_{(j, i)}) \odif{N}_{(j, i)}, \label{eqn:leaky_state} \\
        \vec{y}_{(j, i)}        & = \vec{z}_{(j, i)}, \label{eqn:leaky_output}
      \end{align}
    \end{subequations}
    where \(\vec{z}_{(j, i)} \in \mathbb{R}^{d_{j}}\).
    % A jump resets \(\vec{z}_{(i, j)}\) to \(\vec{x}_j(t^-)\), but leaks at (potentially destabilizing) rate \(a \in \R\), modelling data loss or error accumulations.
    Here, \( a_{(j,i)} \in \mathbb{R} \) is the decay (or potentially amplification, if \( a_{(j,i)} > 0 \)) rate.
    \Cref{fig:example-trajectory-comparison:b} shows how the channel state~\( \vec{z}_{(j,i)} \)  jumps to match \( \vec{x}_j \)  but then decays according to~\( e^{at} \) until the next communication event.
  \end{example}

  Naturally, a multitude of other forms of channel dynamics are possible, including stochastic delay channels, synaptic integration channels, consensus channels or observer channels, but we believe the two exemplified versions are the most common and serve as a good presentational foundation for the remainder of the paper.

  These examples demonstrate how the SDE framework allows us to model the interplay between continuous node computations and discrete communication events, capturing channel imperfections like information decay or even instability.
  The key question, addressed next, is how algorithms designed using this communication model maintain stability and performance.

\subsection{Distributed Gradient Descent}\label{sec:distributed-gradient-descent}
  This section applies the communication model introduced in \Cref{sec:channel-dynamics} to a well-known problem in optimization, namely distributed gradient descent for qudratic cost functions.
  In doing so, we start with a formal definition of the problem, then motivate the construction of distributed optimization dynamics governed by our novel Poisson jump communication channel model and conclude with   a presentation of the main theoretical contributions, each accompanied by discussions concerning their broader implications.

  % define the optimization problem in question
  Our exploration begins with the formal definition of the reference problem, which serves as the foundational example for illustrating our proposed model.
  % Definition of the global quadratic optimization problem.
% This will be included in the main section where the system model is introduced.

\begin{definition}[Quadratic Optimization Problem]\label{def:quadratic-optimization-problem}
  Consider the nominal unconstrained quadratic optimization problem
  \begin{equation}
    \label{eqn:quadratic-problem}
    \min_{\vec{y} \in \R^d} J(\vec{y}) := \tfrac{1}{2} \vec{y}^\transpose \mat{Q} \vec{y} + \vec{q}^\transpose \vec{y} ,
  \end{equation}
  where \({\mat{Q} = \mat{Q}^\transpose \succ 0}\) is symmetric and positive definite and~\({\vec{q} \in \R^d}\).
  The unique optimal solution is given by~\({\vec{y}^* = -\mat{Q}^{-1}\vec{q}}\).
\end{definition}

  % define the nominal gradient descent flow (for reference, discussion, comparison)
  A well-known method for asymptotically solving this quadratic optimization problem is to formulate it as a continuous-time dynamic system, in particular as the gradient descent flow.
  % Definition of the nominal gradient descent flow.
% Used for context and comparison before introducing the distributed stochastic version.

\begin{definition}[Gradient Descent Flow]\label{def:nominal-gradient-descent-flow}
  Associated with~\eqref{eqn:quadratic-problem}, the gradient descent flow is
  % ----------------------------------------------------
  % Beginn der manuellen Korrektur für die Gleichung:
  % 1. Automatische QED-Platzierung von thmtools abschalten:
  % \let\thmendsymbol\relax
  % ----------------------------------------------------
  \begin{equation}
    \label{eqn:nominal-gd}
    \odif{\vec{y}} = -\parens{\mat{Q} \vec{y} + \vec{q}} \odif{t}.
    % \tag{\hflqqed\theequation} %
    \qedhere
  \end{equation}%
\end{definition}
  In anticipation of the parallels to be drawn with a stochastic setup, we have elected to present the gradient descent flow in a differential formulation.

  In view of our desire to split this nominal problem into smaller sub-problems, we can partition \(\vec{y}\) into \(n\) sub-vectors~\(\vec{y}_p\), and partition \(\mat{Q}\) and \(\vec{q}\) accordingly, resulting in
  \begin{align*}
    \mat{Q} & = \begin{pmatrix}
                  \mat{Q}_{11} & \dots  & \mat{Q}_{1n} \\
                  \vdots       & \ddots & \vdots       \\
                  \mat{Q}_{n1} & \dots  & \mat{Q}_{nn}
                \end{pmatrix}, &
    \vec{q} & = \begin{pmatrix}
                  \vec{q}_1 \\
                  \vdots    \\
                  \vec{q}_n
                \end{pmatrix}.
  \end{align*}

  The resulting gradient of the global objective function~\(J\) with respect to an individual sub-vector~\(\vec{y}_i\) can then be precisely expressed as
  \begin{equation}
    \nabla_{\vec{y}_i} J(\vec{y}) = \mat{Q}_{ii} \vec{y}_i + \sum_{j \neq i} \mat{Q}_{ij} \vec{y}_j + \vec{q}_i.
  \end{equation}
  Consequently, the nominal gradient descent flow~\eqref{eqn:nominal-gd}, when considered from the vantage point of each partitioned component~\(\vec{y}_i\), would ideally evolve according to
  \begin{equation}
    \odif{\vec{y}}_i = -\parens{\mat{Q}_{ii} \vec{y}_i + \sum_{j \neq i} \mat{Q}_{ij} \vec{y}_j + \vec{q}_i} \odif{t}.
  \end{equation}
  This formulation can already be thought of as a distributed algorithm, however at each time instant the vectorfield describing the evolution of \(\vec{y}_i\) requires knowledge of the neighbour  sub-vectors~\(\vec{y}_j\).
  In a distributed system, for an agent~\(\nu_i\), the sub-vectors~\(\vec{y}_j\) for \(j \neq i\) will generally be not available, so our idea is to replace them with the channel output~\(\vec{z}_{(j, i)}\) of our proposed channel model instead.
  Jointly with the channel dynamics, this gives rise to our proposed distributed gradient descent flow with Poisson-distributed communication events.

  % define the distributed stochastic gradient descent
  % Definition of the distributed stochastic gradient descent flow with Poisson jumps and drift.
% This is the core system model in the paper.

\begin{definition}%[Distributed Gradient Descent Flow with Poisson-Distributed Communication Events]
  \label{def:distributed-gd-flow}
  Associated with~\eqref{eqn:quadratic-problem},
  the distributed gradient descent flow with Poisson-distributed communication events for agent~\(\nu_i\) with \( i \in \interval{1}{n} \cap \N\) and \(j \in \interval{0}{n}\cap{\N} \setminus \{i\} \), is
  \begin{subequations}
    \label{eqn:def-distributed-gd-distributed-dynamics}%
    \begin{align}
      \label{eqn:def-distributed-gd-primal-dynamics}
      \odif{\vec{x}_i}(t)
       & =
      -\parens{
          \mat{Q}_{ii} \vec{x}_i(t)
          + \sum_{j \neq i} \mat{Q}_{ij} \vec{z}_{(j, i)}(t) + \vec{q}_i
        } \odif{t},
      \\
      \label{eqn:def-distributed-gd-copy-dynamics}
      \odif{\vec{z}_{(j, i)}}(t)
       & =
      \begin{aligned}[t]
          & a_{(j,i)}(t) \vec{z}_{(j, i)}(t)  \odif{t}
        \\
        + & \parens*{
            \vec{x}_j(t) - \vec{z}_{(j, i)}(t)
          } \odif{N}_{(j,i)}(t),
      \end{aligned}
    \end{align}
  \end{subequations}
  where \(\mat{Q}_{ij} \in \R^{d_i \times d_j}\), \(\vec{q}_i \in \R^{d_i}\), \(a_{(j,i)}(t) \in \R \), and \(N_{(j,i)}\) are independent Poisson processes with rates \(\lambda_{(j, i)}\).
\end{definition}

  System~\eqref{eqn:def-distributed-gd-distributed-dynamics} highlights the application of our modelling technique.
  Firstly, channel dynamics~\eqref{eqn:def-distributed-gd-copy-dynamics} for communication events from agents~{\(\nu_j\)} are augmented to each agent~\(\nu_i\) in an existing, well known algorithm, and then the algorithm is adapted to use the local copies from the channel dynamics where applicable~\eqref{eqn:def-distributed-gd-primal-dynamics}.

  For easier analysis, we switch to an equivalent representation using error dynamics relative to the unknown optimal solution~\(\vec{y}^*\).

  %  define the distributed error dynamics
  \begin{lemma}[Error Dynamics]\label{lemma:distributed-error-dynamics}
  % ----------------------------------------------------
  % Beginn der Korrektur
  % \let\oldeqd\qed % Speichern des originalen \qed
  % % Neu definieren: Füge negativen Abstand ein, dann setze \qed (mit \oldeqd)
  % \renewcommand{\qed}{%
  %     \vspace{-2.0\baselineskip}% <--- PASSE DIESEN WERT AN
  %     \oldeqd
  %   }
  \setlength{\belowdisplayskip}{-1.0\baselineskip}% <--- PASSE DIESEN WERT AN
  % ----------------------------------------------------
  Consider the distributed gradient descent flow~\eqref{eqn:def-distributed-gd-distributed-dynamics}, associated with the quadratic optimization problem~\eqref{eqn:quadratic-problem}.
  Let the error be~\(\tilde{\vec{x}} := \vec{x} - \vec{y}^* \in \R^d\).
  Furthermore, let the copy errors be \(\tilde{\vec{z}}_{(j, i)} := \vec{z}_{(j, i)} - \vec{y}^*_j \in \R^{d_j}\) for agent~\(\nu_i\)'s copy of agent \(\nu_j\)'s state~\(\vec{x}_j\), with \( i \in \interval{1}{n} \cap \N\) and \(j \in \interval{0}{n}\cap{\N} \setminus \{i\} \).
  The resultant error dynamics are
  \begin{subequations}%
    \label{eqn:lem:distributed-error-dynamics:errors}
    \begin{align}
      \label{eqn:lem:distributed-error-dynamics:primal-error}
      \odif{\tilde{\vec{x}}}_i(t)
       & =
      -\Big(
      \mat{Q}_{ii} \tilde{\vec{x}}_i(t)
      + \sum_{j \neq i} \mat{Q}_{ij} \tilde{\vec{z}}_{(j, i)}(t)
      \Big) \odif{t},
      \\
      \label{eqn:lem:distributed-error-dynamics:copy-error}
      \odif{\tilde{\vec{z}}}_{(j, i)}(t)
       & =
      \begin{aligned}[t]
          & \Big(
        a_{(j, i)}(t) \tilde{\vec{z}}_{(j, i)}(t)
        + a_{(j, i)}(t) \vec{y}^*_j
        \Big) \odif{t}
        \\
        + & \parens*{\tilde{\vec{x}}_j(t) - \tilde{\vec{z}}_{(j, i)}(t) }\odif{N}_{(i,p)}
      \end{aligned}
    \end{align}
  \end{subequations}
\end{lemma}

  The error system will be used for analysis only, implementation must be done using~\eqref{eqn:def-distributed-gd-distributed-dynamics}, which we will refer to for the upcoming  stability result.
  In particular, we will be considering the stability of a state \({(\vec{x}^*, \vec{z}^*) = (\vec{y}^*, \operatorname{vec}(\vec{y}^*_{-1}, \dots, \vec{y}^*_{-n}))}\), as~\eqref{eqn:def-distributed-gd-distributed-dynamics} does not exhibit a trivial steady-state solution in general.
  In the language of \Cref{def:gpues-ms}, we can state the following theorem.
  % TODO: Raik: "in the language of Definition 4 ...

  % theorem for existence of stabilizing rates
  % Statement for the theorem on Practical Exponential Stability in the Mean-Square Sense (PES-MS).
% This theorem provides sufficient conditions for PES-MS, highlighting the dependency
% of the ultimate bound on the optimal solution y^* and a tunable parameter \rho.

\begin{theorem}
  %[Practical Exponential Stability in the Mean-Square Sense]
  \label{thm:practical-exponential-stability-ms}
  Consider the distributed gradient descent system~\eqref{eqn:def-distributed-gd-distributed-dynamics}.
  Assume that for each \(a_{(j,i)}(t)\), there exists a bound~\({a}_{ji}\), such that \(\abs{a_{(j,i)}(t)} \leq {a}_{ji}\) for all \(t \in \interval[open right]{t_0}{\infty}\).
% 
% TODO Raik: Make it dependent on \gamma?
  % Then, for every \(\parens{\rho_i}_{i=1, \ldots, n}  \in \interval[open]{0}{\infty}^n\),
  Then, 
  \begin{enumerate}[label={\small(\alph*)}]
    \item
        \label{item:thm:practical-exponential-stability-ms:a}
          % stability achievable
          for every \( \gamma \in \interval[open]{0}{\infty}\), there exists a
          constant communication rate \(\lambda_{\text{s}} \in \interval[open right]{0}{\infty}\), such that for all \(\lambda_{(j, i)} \in \interval[open]{\lambda_{\text{s}}}{\infty}\), the state~\((\vec{x}^*, \vec{z}^*)\) is \emph{globally practically uniformly exponentially  stable in the mean-square sense} for~\eqref{eqn:def-distributed-gd-distributed-dynamics} with ultimate bound~\(\gamma\);
    \item
          % convergence rate achievable
          for every \( \gamma \in \interval[open]{0}{\infty}\) and \( \beta \in \interval[open]{0}{2\eigmin{\mat{Q}}} \), there exists a
          constant communication rate \(\lambda_{\text{d}} \in \interval[open right]{\lambda_{\text{s}}}{\infty}\) such that for all choices \(\lambda_{ (j, i)} \in \interval[open]{ \lambda_{\text{d}}}{\infty}\), 
          the stability properties from~\ref{item:thm:practical-exponential-stability-ms:a} hold and the exponential convergence rate in~\eqref{eqn:gpues-ms-lyapunov-condition} is~\(\beta\).
          \qedhere
  \end{enumerate}
\end{theorem}

  \begin{remark}
    \label{rem:pes-ms:1}
    The practical stability bound \( \gamma \) in \Cref{def:gpues-ms} is achieved by selecting positive weighting parameters $(\rho_j)_{j=1..n}$ in the proof (see Step 2).
    These parameters must be chosen to satisfy \( \gamma' \le c_3 \gamma\), where \( c_3 \) is the resulting convergence rate (e.g., \( c_3=\beta \) for statement (b) of~\Cref{thm:practical-exponential-stability-ms}) and \( \gamma' = \sum_{i=1}^{n}\sum_{j\ne i}\frac{1}{\rho_{j}}\norm{ \vec{y}_{j}^{*} }^{2} \).
    This reveals a fundamental trade-off: Achieving a smaller (more accurate) ultimate bound \( \gamma \) requires choosing larger weighting parameters \( \rho_j \).
    Larger $\rho_j$ generally lead to more restrictive bounds in the stability analysis \inversion{full}{(Step 5)}, which in turn demand higher sufficient communication rates \( \lambda_{\text{s}} \) and \( \lambda_{\text{d}} \).
    The size of the ultimate bound, characterized by \( \gamma \), also scales with the number~\(n\) of interfering channels  and the magnitude of the optimal solution \( \vec{y}^* \).
\end{remark}

\begin{remark}~~
    {}{}
    Sufficient choices for the stabilizing rate \( \lambda_{\text{s}} \)  and the convergence-guaranteeing rate \(\lambda_{\text{d}}\) can be analytically obtained. 
    Their derivation, which depends on bounds of the system matrices and the chosen weighting parameters \( (\rho_j)_{j=1, \ldots, n} \), is detailed in the proof of this theorem, see~\eqref{eqn:proof:thm:pes-ms:rates}.
\end{remark}

  The proof of this theorem is deferred to the appendix in favour of a more streamlined presentation.
  Observe, that the theorem states the convergence rate and practical stability  for the variance of the solutions to the error system.
  However, second moment stability entails first moment stability \cite{mao2007stochastic}, i.\,e.\ we also have
  \begin{equation*}
    \E{\norm{\vec{s}}} \leq \sqrt{c_3} \E{\norm*{ \vec{s}^0 }} e^{- \frac 1 2  \beta\parens{t- t_0}} + \sqrt{\gamma}
  \end{equation*}
  by application of Jensen's inequality and by subadditivity of~\(\sqrt{\cdot}\).

  The first part of the theorem states sufficient conditions on the communication channels for achieving stability.
  The second part notably relates the exponential convergence rate to the convergence rate of the nominal problem and it can be observed, that for sufficiently high communication rates the convergence in expectation can be designed to achieve any rate up to the nominal convergence rate from the nominal gradient descent flow in~\Cref{def:nominal-gradient-descent-flow}.
  The overall convergence rates will be limited by \(\eigmin{\mat{Q}}\), for the nominal as well as for the distributed gradient descent.

  % The last part addresses the ball of practical stability.
  % In the representation we found, it is obvious that this region is scaling with the unknown optimal value \(\vec{y}^*\), which makes it a multiplicative disturbance.
  % % TODO: Think about the scaling with M^{-1} and (n-1).
  % Depending on future use-cases, this might be an advantage or an obstacle in comparison to an additive disturbance.
  % Regardless, by selecting scaling factors~\(\rho_i\) per partition, this multiplicative disturbance can be scaled down arbitrarily  per partition~\(i\), however at the cost of driving the bounds on the communication rates up.

  If the communication channel lack drifts, however, we can conveniently state an even tighter result.
  % Corollary to the Existence of Analytically Computable Stabilizing Communication Rates Theorem.
% This corollary states Mean-Square Exponential Stability (MS-ES) when drift rates are zero.

\begin{corollary}%[Mean-Square Exponential Stability]
  \label{cor:mean-square-exponential-stability}
  Consider the distributed gradient descent system~\eqref{eqn:def-distributed-gd-distributed-dynamics} with zero channel drift rates (i.\,e., \(a_{(j, i)}(t) \equiv 0\) for all communication channels~\((j, i)\)).
  Then, in \Cref{thm:practical-exponential-stability-ms} the qualifier \emph{practical} may be removed, i.\,e.~\(\gamma= 0\).
\end{corollary}

  Unfortunately, the sufficient conditions on \(\lambda_{\text{s}}\) and \(\lambda_{\text{d}}\) are only computable by inverting \(\mat{Q}\), solving  large Linear Matrix Inequalities or assuming a particular beneficial property like diagonal dominance of~\(\mat{Q}\), which are an obstacle in a completely distributed setup, where ideally each agent should be able to determine its own communication rates independently from its neighbours.
  As such, our result  remains mostly an existence result, but might be usable nonetheless, as manual tuning of rates or adaptation guarantees stability for sufficiently high rates.

\section{Simulation}\label{sec:simulation}
 
To illustrate our theoretical findings from the previous section, we apply~\eqref{eqn:def-distributed-gd-distributed-dynamics} to an academic example.
We consider a quadratic optimization problem of form~\eqref{eqn:quadratic-problem} with
\begin{align*}
  \mat{Q} & = \begin{bmatrix}
                4 & 2 & 1 & 1 & 0 & 2 \\
                2 & 5 & 2 & 0 & 2 & 1 \\
                1 & 2 & 6 & 3 & 1 & 0 \\
                1 & 0 & 3 & 4 & 2 & 1 \\
                0 & 2 & 1 & 2 & 5 & 2 \\
                2 & 1 & 0 & 1 & 2 & 4
              \end{bmatrix}
          & \text{ and }
          &
          &
  \vec{q} & = \begin{bmatrix}
                -9  \\
                -15 \\
                -22 \\
                -12 \\
                -10 \\
                -5
              \end{bmatrix}
  % &
  % c       & =3
\end{align*}
and formulate our proposed distributed gradient descent flow with Poisson-distributed communication rates for \(n=3\) heterogeneous agents with respective state dimensions \(d_1 = 3\), \(d_2 = 2\), \(d_3 = 1\).
We deliberately have chosen a symmetric positive definite matrix~\(\mat{Q}\), which is not diagonally dominant to highlight this is not a required property.
The nominal convergence rate is characterized by~\(\eigmin{\mat{Q}} \approx \num{0.31}\).

We simulate our proposed algorithm in different scenarios for different choices of communication rates as well as for channels with and without drift.
For the continuous part of the SDE we used an \emph{ode1}-solver with step-size \(h = \num{0.01}\).
Throughout the simulations, we have opted to select all communication and drift rates~\(\lambda=\lambda_i\) equal for a  more compact presentation in the paper.

\begin{figure}
  \centering
  \input{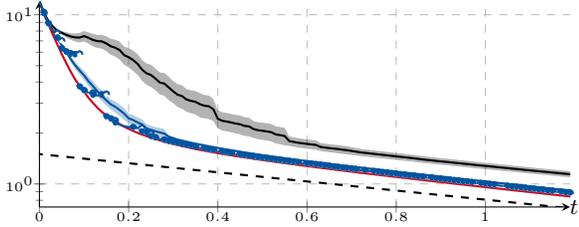}
  \caption{Comparison of Lyapunov-function evolutions with zero channel drift rates~(\(\mat{a} = \mat{0}\)). Depicted are \(\bar{V}(t)\) for diffent communication rates: (\(\lambda_1 = 10\),\ref{plot:figure:lyapunov-function-evolution:A:mean}), (\(\lambda_2 = 27\),\ref{plot:figure:lyapunov-function-evolution:B:mean}) and (\(\lambda_3 = 50\),\ref{plot:figure:lyapunov-function-evolution:C:mean}) alongside their respective confidence intervals, as well as a reference exponential function with decay rate (\(2\eigmin{Q}\),\ref{plot:figure:lyapunov-function-evolution:reference}) and a single sample path~(\ref{plot:figure:lyapunov-function-evolution:B:V50}).
    }\label{fig:lyapunov-function-trajectories}
\end{figure}

\begin{experiment}
  Let~\(\mat{a} = \mat{0}\).
  In view of \Cref{thm:practical-exponential-stability-ms}, we obtain numerically the sufficient communication rate~\(\lambda_{\text{s}} \approx \num{26.56}\).
  Let \(\lambda_1 = 10\), \(\lambda_2 = 27\) and \(\lambda_3 = 50\).
\end{experiment}
The simulation results are depicted in \Cref{fig:lyapunov-function-trajectories}, which shows a comparison of the evolution of the Lyapunov-function~\eqref{eqn:proof-V-candidate} for different communication rates alongside the the nominal convergence bound~\(\frac 3 2 e^{-2 \eigmin{Q} t}\).
A single selected sample path \(V(\vec{s}(\cdot, \omega_1))\) is shown, which showcases, that individual sample paths of \(V\) need not be decreasing at all times in this stochastic setting, unlike for a deterministic system.

Taking the average \(\bar{V}(t) := \frac 1 N \sum_{k=1}^N V(\vec{s}(t, \omega_k))\) over \(N=100\) sample paths of \(V(\vec{s}(t, \omega_k))\) however illustrates the statement of \Cref{thm:practical-exponential-stability-ms} and highlights the exponential stability in the mean-square sense.
A {95-percentile} confidence interval for~\(\bar{V}\) is depicted alongside, to illustrate the spread of trajectories and the exponential decline of the sample path variance.

Regarding the communication rates \(\lambda\), the following can be observed:
\begin{description}
  \item[\(\lambda = \lambda_1 < \lambda_{\text{s}}\)]
        Even with communication rates smaller than \(\lambda_{\text{s}}\), we observe convergence at ultimately nearly identical rates to the optimal solution in some cases, which was not part of the statement in~\Cref{thm:practical-exponential-stability-ms}.
  \item[\(\lambda \in \{\lambda_2, \lambda_3 \}, \lambda_{\text{r}} < \lambda_3  < \lambda_4 \)]
        We observe, that increasing the communication rates beyond \(\lambda_{\text{s}}\) only yields insignificant improvements for the convergence rate of the Lyapunov function, unlike discussed.
\end{description}

\begin{figure}
  \centering
  % --- THE OPTIMIZATION: Read the table ONCE ---
% Read the data from 'simulation_data.csv' and store it in a macro named '\data'
\pgfplotstableread[col sep=comma]{src/figures/figure.lyapunov-function-evolution/l26a1/dgd_stats_data.csv}\dataA
\pgfplotstableread[col sep=comma]{src/figures/figure.lyapunov-function-evolution/l51a1/dgd_stats_data.csv}\dataB

% \pgfplotstableread[col sep=comma]{src/figures/figure.lyapunov-function-evolution/l10a0/dgd_sample_path_run_1.csv}\dataSPA
% \pgfplotstableread[col sep=comma]{src/figures/figure.lyapunov-function-evolution/dgd_sample_path_run_50.csv}\dataSPB
% \pgfplotstableread[col sep=comma]{src/figures/figure.lyapunov-function-evolution/dgd_sample_path_run_100.csv}\dataSPC
% ---------------------------------------------

\begin{tikzpicture}[
    % trim axis left,
  ]
  \begin{axis}[
      xlabel={\(t\)},
      ymode = log,
      grid=major,
      xmin = 0,
      xmax = 29.999,
      %   ymin=-3, ymax=3,
      %   legend pos=south east,
      width = \linewidth,
      height = 0.5*\linewidth,
      % scale only axis,
      axis lines=left,
      ylabel near ticks,
      % ylabel style={at={(ticklabel cs:1)},anchor=south,rotate=-90}, % Move y-label to the top and prevent rotation
      tick label style={font=\tiny},
      yticklabel style={
          inner sep = 0pt,
        },
      label style={font=\small},
    ]

    % % -----------------------------------------------------
    % % 1.1. DEFINE AND NAME THE UPPER PATH (pi_upper)
    % % We use 'draw=none' and 'forget plot' since we only want to use this path for filling.
    % \addplot[
    %     draw=none, % Don't draw the line itself
    %     name path=upper_pi, % <--- NAME THE PATH
    %     forget plot, % Don't include this path in the legend
    %   ] table[
    %     x=time,
    %     y=pi_upper,
    %     col sep=comma,
    %   ] {\dataA};

    % % -----------------------------------------------------
    % % 1.2. DEFINE AND NAME THE LOWER PATH (pi_lower)
    % \addplot[
    %     draw=none,
    %     name path=lower_pi, % <--- NAME THE PATH
    %     forget plot,
    %   ] table[
    %     x=time,
    %     y=pi_lower,
    %     col sep=comma,
    %   ] {\dataA};

    % % -----------------------------------------------------
    % % 1.3. SHADE THE REGION BETWEEN THE NAMED PATHS
    % \addplot[
    %     fill=black, % Color of the fill
    %     opacity=0.3, % Make it semi-transparent
    %   ] fill between[
    %     of=lower_pi and upper_pi, % <--- FILL BETWEEN THE TWO NAMED PATHS
    %   ];
    % -----------------------------------------------------
    % 2.1. DEFINE AND NAME THE UPPER PATH (ci_upper)
    % We use 'draw=none' and 'forget plot' since we only want to use this path for filling.
    \addplot[
        draw=none, % Don't draw the line itself
        name path=upper_ci, % <--- NAME THE PATH
        forget plot, % Don't include this path in the legend
      ] table[
        x=time,
        y=ci_upper,
        col sep=comma,
      ] {\dataA};
    % -----------------------------------------------------
    % 2.2. DEFINE AND NAME THE LOWER PATH (ci_lower)
    \addplot[
        draw=none,
        name path=lower_ci, % <--- NAME THE PATH
        forget plot,
      ] table[
        x=time,
        y=ci_lower,
        col sep=comma,
      ] {\dataA};
    % -----------------------------------------------------
    % 2.3. SHADE THE REGION BETWEEN THE NAMED PATHS
    \addplot[
        fill=black, % Color of the fill
        opacity=0.3, % Make it semi-transparent
      ] fill between[
        of=lower_ci and upper_ci, % <--- FILL BETWEEN THE TWO NAMED PATHS
      ];\label{plot:figure:lyapunov-function-evolution-drift:A:ci}
    % -----------------------------------------------------
    % 3. PLOT THE MEAN LINE (last so it appears on top)
    \addplot [
        black,
        thick,
      ] table [
        % skip first n = 13,
        x = time,
        y = mean,
        col sep = comma,
      ] {\dataA};\label{plot:figure:lyapunov-function-evolution-drift:A:mean}

    % % -----------------------------------------------------
    % % 1.1. DEFINE AND NAME THE UPPER PATH (pi_upper)
    % % We use 'draw=none' and 'forget plot' since we only want to use this path for filling.
    % \addplot[
    %     draw=none, % Don't draw the line itself
    %     name path=upper_pi, % <--- NAME THE PATH
    %     forget plot, % Don't include this path in the legend
    %   ] table[
    %     x=time,
    %     y=pi_upper,
    %     col sep=comma,
    %   ] {\dataB};

    % % -----------------------------------------------------
    % % 1.2. DEFINE AND NAME THE LOWER PATH (pi_lower)
    % \addplot[
    %     draw=none,
    %     name path=lower_pi, % <--- NAME THE PATH
    %     forget plot,
    %   ] table[
    %     x=time,
    %     y=pi_lower,
    %     col sep=comma,
    %   ] {\dataB};
    % % -----------------------------------------------------
    % % 1.3. SHADE THE REGION BETWEEN THE NAMED PATHS
    % \addplot[
    %     fill=blue, % Color of the fill
    %     opacity=0.3, % Make it semi-transparent
    %   ] fill between[
    %     of=lower_pi and upper_pi, % <--- FILL BETWEEN THE TWO NAMED PATHS
    %   ];
    % -----------------------------------------------------
    % 2.1. DEFINE AND NAME THE UPPER PATH (ci_upper)
    % We use 'draw=none' and 'forget plot' since we only want to use this path for filling.
    \addplot[
        draw=none, % Don't draw the line itself
        name path=upper_ci, % <--- NAME THE PATH
        forget plot, % Don't include this path in the legend
      ] table[
        x=time,
        y=ci_upper,
        col sep=comma,
      ] {\dataB};
    % -----------------------------------------------------
    % 2.2. DEFINE AND NAME THE LOWER PATH (ci_lower)
    \addplot[
        draw=none,
        name path=lower_ci, % <--- NAME THE PATH
        forget plot,
      ] table[
        x=time,
        y=ci_lower,
        col sep=comma,
      ] {\dataB};
    % -----------------------------------------------------
    % 2.3. SHADE THE REGION BETWEEN THE NAMED PATHS
    \addplot[
        fill=blue, % Color of the fill
        opacity=0.3, % Make it semi-transparent
      ] fill between[
        of=lower_ci and upper_ci, % <--- FILL BETWEEN THE TWO NAMED PATHS
      ];\label{plot:figure:lyapunov-function-evolution-drift:B:ci}
    % -----------------------------------------------------
    % 3. PLOT THE MEAN LINE (last so it appears on top)
    \addplot [
        blue,
        thick,
      ] table [
        % skip first n = 13,
        x = time,
        y = mean,
        col sep = comma,
      ] {\dataB};\label{plot:figure:lyapunov-function-evolution-drift:B:mean}

  \end{axis}
\end{tikzpicture}
  \caption{Comparison of Lyapunov-function evolutions with unstable channels~(\(\mat{a}_{ji} = 1\) for \(j \neq i\)). Depicted are \(\bar{V}(t)\) for different communication rates (\(\lambda_1 = 26\),\ref{plot:figure:lyapunov-function-evolution-drift:A:mean}), (\(\lambda_2 = 51\),\ref{plot:figure:lyapunov-function-evolution-drift:B:mean}).}\label{fig:lyapunov-function-trajectories-drift}
\end{figure}
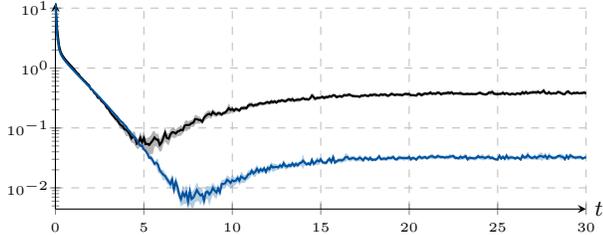

\begin{experiment}
  Let~\(\mat{a}_{ip} = 1\) for each \(i \neq p\) and \(\mat{a}_{ii} = 0\).
  % Let  \(\rho_p = 1\) for all \(p \in \interval{1}{n} \cap \N\).
  In view of \Cref{thm:practical-exponential-stability-ms}, we obtain numerically the sufficient communication rate~\(\lambda_{\text{s}} \approx \num{50.57}\).
  Let \(\lambda_1 = 26\) and \(\lambda_2 = 51\).
\end{experiment}
The simulation results are depicted in \Cref{fig:lyapunov-function-trajectories-drift}, which shows a comparison of the evolution of the Lyapunov-function~\eqref{eqn:proof-V-candidate} for different communication rates.
Once again, the averages \(\bar{V}(t)\) and {95-percentile} confidence interval for~\(\bar{V}(t)\) are depicted.
Regarding the communication rates \(\lambda\), the following can be observed:
\begin{description}
  \item[\(\lambda = \lambda_1 < \lambda_{\text{s}}\)]
        Even with communication rates smaller than \(\lambda_{\text{s}}\), we observe practical convergence, which was not part of the statement in ~\Cref{thm:practical-exponential-stability-ms}.
  \item[\(\lambda = \lambda_2 > \lambda_{\text{s}} \)]
        Increasing the communication rates  means essentially selecting smaller ultimate bound~\(\gamma\) (or equivalently: higher scaling weights \(\rho\)) and renders the disk of practical convergence smaller, as stated by~\Cref{rem:pes-ms:1}.
\end{description}

\section{Discussion, Conclusion \& Outlook}\label{sec:discussion}
 In this paper, we have presented a novel framework for designing communication-aware distributed optimization algorithms.
The core of our approach, motivated by networked systems and modern hardware, is to model the communication channel itself as a dynamic, stochastic process using SDEs driven by Poisson Jumps.
This fundamentally shifts the design paradigm: instead of treating communication scarcity as an afterthought imposed on an existing algorithm—as is common in many event-triggered schemes—our framework embeds these channel constraints into the mathematical flow from its inception.
This allows us to take nominal continuous-time dynamics (representing ideal, fully-connected systems) and rigorously analyze their stability and performance under sparse, asynchronous, and stochastic communication events, such as exponentially distributed sampling times.
Unlike state-dependent event-triggering, this Poisson-driven approach does not require agents to continuously monitor local state changes to make communication decisions, which significantly reduces the computational and energy overhead associated with triggering logic.
We have illustrated the process of turning a continuous-time dynamic system (the gradient descent flow) into a network system with constant communication (just a matter of different interpretation of the gradient descent flow), then applying our proposed  channel dynamics model to obtain the distributed gradient descent flow with Poisson-distributed communication events.
Using our modelling technique, we were then able to state stability guarantees depending on sufficiently high communication rates in two increasingly meaningful ways: 
Firstly, we show existence of sufficient stabilizing communication rates, and secondly, we show the existence of sufficient  communication rates to achieve a desired convergence performance up to the nominal gradient descent flow.
While it serves as an example to illustrate our modelling technique, the proposed distributed gradient descent flow with Poisson-distributed communication events is a usable algorithm on its own.
Even though the chosen illustrative example (the distributed gradient descent flow) requires a complete communication graph due to its reliance on global state aggregation, this is not a limitation of the SDE framework itself.
The framework's core strength lies in its ability to model and analyze complex (deterministic and stochastic) channel dynamics, including time-varying or unstable channels or sparse random communication patterns over directed graphs.
Applying the Poisson-driven SDE approach to consensus-based algorithms designed for sparse topologies (e.\,g., DGD or ADMM) presents one natural direction for future work.

Another future research direction is to apply the approach to algorithms that are natively designed for sparse communication graphs, such as consensus-based DGD or ADMM.\@
This would directly connect our theoretical framework to the sparsely-connected hardware architectures that motivated this work.
Further, the framework's flexibility invites extensions to more complex problem classes, including constrained quadratic programming and the integration of state observers to handle partial information.

Moreover, a promising direction for future work involves exploring accelerated convergence within this distributed SDE framework.
Leveraging techniques from convex synthesis for optimization algorithms for centralised systems~\cite{scherer2021convex}, could enable the co-design of both the algorithm dynamics and the communication strategy.
This might lead to characterising minimal communication rates to achieve performance exceeding the nominal gradient flow, even under sparse, stochastic communication constraints.
Finally, our focus on continuous-time dynamics positions could be transferred to discrete-time models with jump dynamics, such as~\cite{ebenbauer2025piecewise}.
Investigating the connections and potential translations between these continuous and discrete-time stochastic frameworks represents a fertile ground for further future theoretical exploration, promising a deeper understanding of distributed optimization.

 %  \addtolength{\textheight}{-3cm}   % This command serves to balance the column lengths

 %%%%%%%%%%%%%%%%%%%%%%%%%%%%%%%%%%%%%%%%%%%%%%%%%%%%%%%%%%%%%%%%%%%%%%%%%%%%%%%
 % References
 %%%%%%%%%%%%%%%%%%%%%%%%%%%%%%%%%%%%%%%%%%%%%%%%%%%%%%%%%%%%%%%%%%%%%%%%%%%%%%%
 %  \bibliographystyle{ieeetr}
 \bibliographystyle{IEEEtran}
 \bibliography{references}             % bib file to produce the bibliography

 \appendices%

 % Grab forced line break - \\* - and replace with :
 \renewcommand{\thesectiondis}[2]{\Alph{section}:}

% \section{Definitions}\label{sec:definitions}

\section{Lemmata}\label{sec:lemmata}

 \inversion{full}{\begin{lemma}
  \label{lemma:matrix-bounds}
  Let \( M \) be a symmetric block matrix of the form
  \begin{align*}
    \mat{M} &= \begin{bmatrix}
    \mat{M}_{11} & \mat{M}_{12} \\
    \mat{M}_{12}^\transpose & \mat{\Lambda} + \mat{R} 
    \end{bmatrix}
  \end{align*}
  where
  %   \begin{itemize*}[label={}, mode=unboxed, itemjoin={{, }}, itemjoin*={{ }}]
  % \item
  \( \mat{M}_{11} \) is symmetric and positive definite,  \( \mat{R}\) is symmetric 
  % \item
  % \( {\mat{M}_{21} = \mat{M}_{12}^\transpose }\) and
  % \item
  % \({ \mat{M}_{22} = \mat{\Lambda} + \mat{R} }\)
  and \({ \mat{\Lambda} = \text{diag}(\lambda_k) }\).
  %   \end{itemize*}
  % 
  % and \({ \mat{S} := \mat{\Lambda} + \mat{R} - \mat{K} }\) be the Schur complement of~\( \mat{M}_{11} \).
  The following two implications hold:

  \begin{enumerate}[
      label={\tiny (\arabic*)},
      % wide,
      % labelindent=0pt,
      labelsep=0.5em
    ]
    \item 
    Let~\({ \mat{K} := \mat{M}_{12}^\transpose\mat{M}_{11}^{-1}\mat{M}_{12} }\).
    If all \( \lambda_k \) satisfy
          \[
            \lambda_k > \lambda_{\text{s}} :=  -\eigmin{\mat{R} - \mat{K}},
          \]
    then the matrix~\( \mat{M} \) is positive definite.

    \item 
    % {Eigenvalue Saturation}: 
    Let \( {\mu \in \interval[open]{0}{\eigmin{\mat{M}_{11}}}}\) and~\({ \tilde{\mat{K}}_\mu := \mat{M}_{12}^\transpose\parens{\mat{M}_{11} - \mu \mat{I}}^{-1}\mat{M}_{12} }\).
    % and  \({ \tilde{\mat{K}} := \mat{M}_{12}^\transpose\parens{\mat{M}_{11} - \mu \mat{I}}^{-1}\mat{M}_{12} }\).
    If all \(\lambda_k\) satisfy
          \[
            \lambda_k \ge \lambda_{\text{d}} := \mu - \eigmin{\mat{R}-\tilde{\mat{K}}_\mu},
          \]
    % then the minimal eigenvalue~\(\eigmin{\mat{M}}\) of~\(\mat{M}\) is \({ \eigmin{\mat{M}} = \eigmin{\mat{M}_{11}} }\).
    then \(\mat{M} \succeq \mu \mat{I}    \).
    \qedhere
  \end{enumerate}
\end{lemma}}

\section{Proofs}\label{sec:proofs}
 \inversion{full}{\begin{proof}[of \Cref{lemma:matrix-bounds}]
  \leavevmode
  \begin{enumerate}[
      label={\tiny (\arabic*)},
      wide,
      labelindent=0pt,
    ]
    \item\label{item:proof:lemma:matrix-bounds:a}
          A symmetric matrix \( \mat{M} \) is positive definite if and only if its principal submatrix \( \mat{M}_{11} \) is positive definite and its Schur complement \( \mat{S}:= \mat{\Lambda} + \mat{R} - \mat{K}  \) is positive definite~\cite[p.~495]{horn2012matrix}.
          According to the assumptions of this lemma,~\( \mat{M}_{11} \) is positive definite.
          Thus, we need to show that~\( \mat{S} \) is positive definite, which is equivalent to its smallest eigenvalue \( \eigmin{\mat{S}} \) being strictly positive.
          We lower-bound~\(\eigmin{\mat{S}}\) via
          \[
            \eigmin{\mat{S}} = \eigmin{\mat\Lambda + (\mat{R} - \mat{K})} \ge \eigmin{\mat\Lambda} + \eigmin{\mat{R} -\mat{K}}
          \]
          using Weyl's inequality for eigenvalues of sums of symmetric matrices~\cite[Corollary~4.3.15]{horn2012matrix}.

          Since \({ \mat\Lambda = \operatorname{diag}(\lambda_k) }\), its smallest eigenvalue is \( {\eigmin{\mat\Lambda} = \min_k(\lambda_k)} \).
          Therefore, for \( \mat{S} \) to be positive definite, we require
          \[
            \min_k(\lambda_k) + \eigmin{\mat{R} - \mat{K}} > 0.
          \]
          % This condition is satisfied if all \( \lambda_k \) satisfy the condition from this lemma's statement.
          Rearranging yields the first statement of this lemma.

    \item

        The condition \( \mat{M} \succeq \mu \mat{I} \) is equivalent to showing \( \mat{M}' := \mat{M} - \mu \mat{I} \succeq 0\).
        Since \( \mu < \eigmin{\mat{M}_{11}} \), the block \( \mat{M}'_{11} \succ 0\) and subtracting the diagonal matrix \( \mu \mat{I} \) inherits the symmetry from \( \mat{M} \)  to \( \mat{M}' \).
        Applying the first statement of this lemma to \( \mat{M}' \) and considering positive semi-definiteness of \(\mat{S}\) yields the second statement of this lemma.
          % First, we establish bounds for
          % \(
          % \eigmin{\mat{M}} \le \eigmin{\mat{M}_{11}}
          % \)
          % via the properties of the Rayleigh quotient applied to the principal submatrix~\(\mat{M}_{11}\).
          % Next, a Schur complement result for symmetric block matrices states
          % \(
          %   \eigmin{\mat{M}} \ge \min\{\eigmin{\mat{M}_{11}}, \eigmin{\mat{S}}\}.
          % \)
          % From part~\ref{item:proof:lemma:matrix-bounds:a},  \( \eigmin{\mat{S}} \ge \min_k(\lambda_k) + \eigmin{\mat{R} - \mat{K}} \).
          % Now, consider the condition \({\min_k (\lambda_k) \ge \lambda_{\text{r}}} \) from the statement of this lemma.
          % If this condition holds, then,
          % \[
          %   \min_k(\lambda_k) + \eigmin{\mat{R}-\mat{K}} \ge \eigmin{\mat{M}_{11}},
          % \]
          % which implies \( \eigmin{\mat{S}} \ge \eigmin{\mat{M}_{11}} \).
          % Substituting this into the lower bound for \( \eigmin{\mat{M}} \) we obtain
          % \(
          %   \eigmin{\mat{M}} \geq  \eigmin{\mat{M}_{11}}
          % \)
          % which concludes the second statement of this lemma.
          \qedhere
          %   For the backward direction, assume \( \eigmin{\mat{M}} = \eigmin{\mat{M}_{11}} \).
          %   If \( \min_k \lambda_k < \lambda_{\text{r}} \), then \( \eigmin{\mat{S}} < \eigmin{\mat{M}_{11}} \).
          %   In this case, the eigenvector corresponding to \( \eigmin(M) \) may not be exclusively in the \( M_{11} \) block, and the coupling terms \( M_{12} \) and \( M_{21} \) would typically cause \( \eigmin(M) \) to be strictly less than \( \eigmin(M_{11}) \). Thus, for equality to hold, the condition \( \min_k \lambda_k \ge \lambda_{recovery} \) must be satisfied.
  \end{enumerate}
\end{proof}}
 \inversion{full}{% Proof for Theorem: Practical Exponential Stability in Mean-Square Sense (\ref{thm:practical-exponential-stability-ms})

\begin{proof}[of \Cref{thm:practical-exponential-stability-ms}]\label{proof:pes-ms}
  This proof aims to show that the system is globally practically uniformly exponentially  stable in the mean-square sense by applying the Lyapunov theorem from \Cref{lemma:gpues-ms-lyapunov-like-theorem}.
  % The strategy involves
  % \begin{enumerate}[
  %     label={\tiny (Step \arabic*)},
  %     wide,
  %     labelindent=0pt,
  %   ]
  %   \item defining  a Lyapunov function candidate,
  %   \item meticulously calculating  the infinitesimal generator,
  %   \item reformatting the infinitesimal generator into a quadratic form,
  %   \item applying \Cref{lemma:matrix-bounds} to show that the conditions of the Lyapunov theorem are met for sufficiently large communication rates, and lastly
  %   \item discussing the  radius of ultimate convergence.
  % \end{enumerate}

Let \( \gamma \in \interval[open]{0}{\infty} \) be given. 
For the second statement, let \( \beta \in \interval[open]{0}{2\eigmin{\mat{Q}}} \) also be given.
  % \noindent
  % Let \(\parens{\rho_i}_{i=1, \ldots, n}  \in \interval[open]{0}{\infty}^n\) arbitrary but fixed and let
  % \begin{equation*}
  %   \mat{\mathrm{P}} := \operatorname{diag}\parens*{\rho_{1}\mat{I}_{d_1}, \dots, \rho_{n}\mat{I}_{d_n}},
  % \end{equation*}
  \begin{enumerate}[
      label={\tiny (Step \arabic*)},
      wide,
      labelindent=0pt,
    ]
    \item
          For notational convenience, define firstly \(
          {
              \vec{s}^\transpose := \begin{bsmallmatrix}
                \tilde{\vec{x}}^\transpose
                &
                \vec{e}^\transpose
              \end{bsmallmatrix}
            }
          \),          where \(\vec{e}\) is the stacked vector of copy errors, with \({\vec{e}_{(j, i)} := \tilde{\vec{x}}_j - \tilde{\vec{z}}_{(j, i)} \in \R^{d_j}}\).
          As a candidate Lyapunov function, we choose the quadratic form
          \begin{equation}
            V(\vec{s}) := \norm{\vec{s}}^2
            = {\sum_{i=1}^n \norm{\tilde{\vec{x}}_i}^2
            + \sum_{i=1}^n \sum_{j \neq i} \norm{\vec{e}_{(j, i)}}^2}.
            \label{eqn:proof-V-candidate}
          \end{equation}
          This function satisfies the first condition of \Cref{lemma:gpues-ms-lyapunov-like-theorem} with \(c_1 = 1\) and \(c_2 = 1\).
          According to \Cref{lemma:gpues-ms-lyapunov-like-theorem}, we need to show that there exist constants \(c_3 > 0\) and \( \gamma' > 0 \) such that \( \mathcal{L}V(\vec{s}) \le -c_3 \norm{\vec{s}}^2 + \gamma' \), where the resulting stability parameters satisfy $\beta= \tfrac{c_3}{c_2} = c_3$ and $\gamma' = {c_3}\gamma$.
    \item
          We directly calculate the infinitesimal generator
          % 2. Compute the infinitesimal generator of the Lyapunov function.
          \begin{align}
             \mathrlap{\mathcal{L}V(\vec{s})} \hphantom{A} & \\
            %  & \labelrel{=}{eqn:infinitesimal-generator:exakt:1}
            % \parens*{\mathcal{L}V}_{\text{(cont)}}
            % + \parens*{\mathcal{L}V}_{\text{(jump)}}
            % \notag
            % \\
             & \labelrel{=}{eqn:infinitesimal-generator:exakt:2}
            \begin{aligned}[t]
               & \sum_{i=1}^n \pdv{V}{\tilde{\vec{x}}_i}^\transpose 
               % \odif{\tilde{\vec{x}}_i}
               \Big( -
      \mat{Q}_{ii} \tilde{\vec{x}}_i(t)
      - \sum_{j \neq i} \mat{Q}_{ij} \tilde{\vec{z}}_{(j, i)}(t)
      \Big)
               \\
              & + \sum_{i=1}^n \sum_{j \neq i} \pdv{V}{\tilde{\vec{z}}_{(j, i)}}^\transpose 
              % \odif{\tilde{\vec{z}}_{(j, i)}}
              \Big(
        a_{(j, i)}(t) \tilde{\vec{z}}_{(j, i)}(t)
        + a_{(j, i)}(t) \vec{y}^*_j
        \Big)
              \\
               & + \sum_{i=1}^n \sum_{j \neq i} \lambda_{(j,i)} \parens*{V(\vec{s} + \Delta\vec{s}) - V(\vec{s})}
            \end{aligned}
            \notag
            \\
             & \labelrel{=}{eqn:infinitesimal-generator:exakt:3}
            \begin{aligned}[t]
               & 2\sum_{i=1}^n {\Bigl(
                \tilde{\vec{x}}_i
                +  \sum_{j \neq i} \parens*{\tilde{\vec{x}}_i - \tilde{\vec{z}}_{(i, j)}}
                \!
                \Bigr)}^\transpose
              \Bigl(
                -\mat{Q}_{ii}\tilde{\vec{x}}_i - \!\!\sum_{j \neq i} \mat{Q}_{ij}\tilde{\vec{z}}_{(j, i)}
                \!
                \Bigr)
              \\
               & -2 \sum_{i=1}^n \sum_{j \neq i}
              \parens*{\tilde{\vec{x}}_j - \tilde{\vec{z}}_{(j, i)}}^\transpose
              \parens*{a_{(j,i)}(t)\tilde{\vec{z}}_{(j, i)} + a_{(j,i)}(t)\vec{y}_j^*}
              \\
               & - \sum_{i=1}^n \sum_{j \neq i} \lambda_{(j,i)}  \parens*{\tilde{\vec{x}}_j - \tilde{\vec{z}}_{(j, i)}}^\transpose \parens*{\tilde{\vec{x}}_j - \tilde{\vec{z}}_{(j, i)}}
            \end{aligned}
            \notag
            \\
             & \labelrel{=}{eqn:infinitesimal-generator:exakt:4}
            \begin{aligned}[t]
               & -2\sum_{i=1}^n \tilde{\vec{x}}_i^\transpose \mat{Q}_{ii}\tilde{\vec{x}}_i
              \\
               & -2\sum_{i=1}^n \sum_{j \neq i} \tilde{\vec{x}}_i^\transpose \mat{Q}_{ij} \parens*{\tilde{\vec{x}}_j - \parens*{\tilde{\vec{x}}_j - \tilde{\vec{z}}_{(j, i)}}}
              \\
               & -2\sum_{i=1}^n \sum_{j \neq i} \parens*{\tilde{\vec{x}}_i - \tilde{\vec{z}}_{(i, j)}}^\transpose \mat{Q}_{ii}\tilde{\vec{x}}_i
              \\
               & -2\sum_{i=1}^n \sum_{j \neq i} \sum_{k \neq i} \parens*{\tilde{\vec{x}}_i - \tilde{\vec{z}}_{(i, j)}}^\transpose \mat{Q}_{ik}\parens*{\tilde{\vec{x}}_k - \parens*{\tilde{\vec{x}}_k - \tilde{\vec{z}}_{(k, i)}}}
              \\
               & -2 \sum_{i=1}^n \sum_{j \neq i}
              \parens*{\tilde{\vec{x}}_j - \tilde{\vec{z}}_{(j,i)}}^\transpose
              a_{(j, i)}(t)
              \parens*{\tilde{\vec{x}}_j -\parens*{\tilde{\vec{x}}_j - \tilde{\vec{z}}_{(j, i)}}}
              \\
               & -2 \sum_{i=1}^n \sum_{j \neq i}
              \parens*{\tilde{\vec{x}}_j - \tilde{\vec{z}}_{(j, i)}}^\transpose
                { a_{(j,i)}(t)\vec{y}_j^*}
              \\
               & - \sum_{i=1}^n \sum_{j \neq i} \lambda_{(j,i)}  \parens*{\tilde{\vec{x}}_j - \tilde{\vec{z}}_{(j, i)}}^\transpose \parens*{\tilde{\vec{x}}_j - \tilde{\vec{z}}_{(j, i)}}
            \end{aligned}
            \notag
            \\
             & \labelrel{=}{eqn:infinitesimal-generator:exakt:5}
            \begin{aligned}[t]
               & -2\sum_{i=1}^n \tilde{\vec{x}}_i^\transpose \mat{Q}_{ii}\tilde{\vec{x}}_i
              - 2\sum_{i=1}^n \sum_{j \neq i} \tilde{\vec{x}}_i^\transpose \mat{Q}_{ij}\tilde{\vec{x}}_j
              \\
               & + 2\sum_{i=1}^n \sum_{j \neq i} \tilde{\vec{x}}_i^\transpose \mat{Q}_{ij}\parens*{\vec{e}_{(j, i)}}
              \\
               & -2\sum_{i=1}^n \sum_{j \neq i} \parens*{ \vec{e}_{(i, j)} }^\transpose \Bigl( \sum_{k=1}^n \mat{Q}_{ik}\tilde{\vec{x}}_k \Bigr)
              \\
               & +2\sum_{i=1}^n \sum_{j \neq i} \sum_{k \neq i} \parens*{\vec{e}_{(i, j)}}^\transpose \mat{Q}_{ik}\parens*{\vec{e}_{(k, i)}}
              \\
               & -2 \sum_{i=1}^n \sum_{j \neq i} \parens*{\vec{e}_{(j, i)}}^\transpose a_{(j,i)}(t)\tilde{\vec{x}}_j
              \\
               & + \sum_{i=1}^n \sum_{j \neq i} \parens*{ 2a_{(j, i)}(t) - \lambda_{(j, i)} } \parens*{\vec{e}_{(j, i)}}^\transpose \parens*{\vec{e}_{(j, i)}}
              \\
               & -2 \sum_{i=1}^n \sum_{j \neq i} \parens*{\vec{e}_{(j, i)}}^\transpose a_{(j, i)}(t)\vec{y}_j^*
            \end{aligned}
            \\
             & \labelrel{=}{eqn:infinitesimal-generator:exakt:7}
            W_1 + 2\,W_2 + W_3 + 2\,W_4 + 2\,W_5 + W_6 + W_7
            \label{eqn:proof_LV_expanded_exact}
          \end{align}
          with expressions
          \begin{subequations}
            \begin{align}
              W_1
                  & = -2 \tilde{\vec{x}}^\transpose \mat{Q}  \tilde{\vec{x}}
              \\
              W_2
                  & = \phantom{+} \sum_{i=1}^n \sum_{j \neq i} \parens*{\vec{e}_{(j, i)}}^\transpose \mat{Q}_{ji}  \tilde{\vec{x}}_i
              \\
              W_3 & = \phantom{+}2\sum_{i=1}^n \sum_{j \neq i} \sum_{k \neq j} \parens*{\vec{e}_{(j, i)}}^\transpose \mat{Q}_{jk}\parens*{\vec{e}_{(k, j)}}
              \\
              W_4
                  & = -\sum_{i=1}^n \sum_{j \neq i} \parens*{ \vec{e}_{(j, i)} }^\transpose \Bigl( \sum_{k=1}^n \mat{Q}_{jk}\tilde{\vec{x}}_k \Bigr)
              \\
              W_5
                  & = - \sum_{i=1}^n \sum_{j \neq i} \parens*{\vec{e}_{(j, i)}}^\transpose a_{(j,i)}(t)\mat{I}_{d_j}\tilde{\vec{x}}_j
              \\
              W_6
                  & =
              \phantom{+}\sum_{i=1}^n \sum_{j \neq i} \parens*{\vec{e}_{(j, i)}}^\transpose \parens*{ 2a_{(j, i)}(t) - \lambda_{(j, i)} } \mat{I}_{d_j} \parens*{\vec{e}_{(j, i)}}
              \\
              \label{eqn:infinitesimal-generator:W7}
              W_7
                  & =
              - 2\sum_{i=1}^n \sum_{j \neq i} \parens*{\vec{e}_{(j, i)}}^\transpose a_{(j, i)}(t)\mat{I}_{d_j}\vec{y}_j^*
            \end{align}
            using relabelled counters where we,
            \begin{steps}
              % \item[\eqref{eqn:infinitesimal-generator:exakt:1}] stress the individual contributions from continuous dynamics and asynchronous jump events,
              \item[\eqref{eqn:infinitesimal-generator:exakt:2}] apply the Itô jump rule, carefully considering the partial derivatives of the Lyapunov function and the value of \(V\) immediately following a jump event and insert~\eqref{eqn:lem:distributed-error-dynamics:errors},
              \item[\eqref{eqn:infinitesimal-generator:exakt:3}] substitute in the error dynamics derived from \Cref{lemma:distributed-error-dynamics},
              \item[\eqref{eqn:infinitesimal-generator:exakt:4}] expand the first term of~\eqref{eqn:infinitesimal-generator:exakt:3} and introduce zero-sum terms to enhance clarity,
              \item[\eqref{eqn:infinitesimal-generator:exakt:5}] expand and combine expressions, explicitly incorporating the newly defined error symbols, such as \(\vec{e}_{(i, p)}\), for improved readability and conciseness and
              \item[\eqref{eqn:infinitesimal-generator:exakt:7}] group the resulting bilinear and quadratic forms, thus expressing the total infinitesimal generator as a summation over individual contributions.
            \end{steps}

            % unneeded, follows in explanation step
            % We apply Young's inequality \(2ab \le \rho a^2 + \frac{1}{\rho} b^2\) to \(2W_7\).
            % We introduce a positive constant \(\rho_{pi}\) (specific to each pair \((p,i)\)) to distribute coefficients effectively.
            Additionally, we upper-bound  the bilinear expression \(W_7\) with the intention to absorb it partly into the bilinear terms as
            \begin{align}
              W_7 & \labelrel{=}{eqn:W7:bound:step:1}
              -2 \sum_{i=1}^n \sum_{j \neq i} \parens*{a_{(j, i)}(t)\mat{I}_{d_j}\vec{e}_{(j, i)}}^\transpose \vec{y}_j^*
              \notag
              \\
                  & \labelrel{\le}{eqn:W7:bound:step:3}
              \sum_{i=1}^n \sum_{j \neq i} \parens*{
                  \rho_{j} \norm{a_{(j,i)}(t)\mat{I}_{d_j}\vec{e}_{(j, i)}}^2
                  + \tfrac{1}{\rho_{j}} \norm{\vec{y}_j^*}^2
                }
              \notag
              \\
              \label{eqn:infinitesimal-genarator:W7:upper-bound}
                  & \labelrel{=}{eqn:W7:bound:step:4}
              \sum_{i=1}^n \sum_{j \neq i}
              \parens*{
                  \rho_{j} a_{(j,i)}(t)^2
                  \norm{\vec{e}_{(j, i)}}^2
                  + \tfrac{1}{\rho_{j}} \norm{\vec{y}_j^*}^2
                }.
            \end{align}
            where we
            \begin{steps}
              \item[\eqref{eqn:W7:bound:step:1}] plug in the original expression for \(W_7\) from~\eqref{eqn:infinitesimal-generator:W7}, with the drift rates~\(a_{(j,i)}\) grouped with the error terms~\(e\),
              \item[\eqref{eqn:W7:bound:step:3}] apply Young's inequality introduce positive weighting constants $(\rho_{j})_{j=1, \ldots, n}$ and
              \item[\eqref{eqn:W7:bound:step:4}] expand the norm using submultiplicativity of norms.
            \end{steps}
          \end{subequations}

          We are now ready to reformulate an upper bound for the infinitesimal generator as a bilinear form.

    \item
          The infinitesimal generator can be upper-bounded by the quadratic form \( -\vec{s}^\transpose \mat{M}(t) \vec{s} + \gamma' \).
          The matrix \(\mat{M}(t)\), encapsulating the quadratic and bilinear terms, takes the block structure
          \begin{align*}
            \mat{M}(t) & = \begin{bmatrix}
                             \mat{M}_{11}    & \mat{M}_{21}^\transpose(t) \\
                             \mat{M}_{21}(t) & \mat{M}_{22}(t)
                           \end{bmatrix}.
          \end{align*}
          The blocks of \(\mat{M}(t)\) are formed by systematically collecting the coefficients from the expanded expression of \(\mathcal{L}V\) from~\eqref{eqn:proof_LV_expanded_exact}, incorporating the upper bound for \(W_7\) from~\eqref{eqn:infinitesimal-genarator:W7:upper-bound}.

          The block \(\mat{M}_{11}\), capturing quadratic terms involving only~\(\tilde{\vec{x}}\), is given by
          \begin{align}
            \mat{M}_{11}
             & = 2\mat{Q}
            \label{eqn:M11_def}
          \end{align}
          The block \(\mat{M}_{21}(t)\) captures bilinear terms involving products of~\(\vec{e}\) and~\(\tilde{\vec{x}}\), defined by contributions from \(W_2\), \(W_4\), and \(W_5\) as
          \begin{align}
            \parens*{\mat{M}_{21}}_{(j, i), (k)}(t)
             & =
            \mat{Q}_{jk}
            \begin{aligned}[t]
               & -
              \begin{cases}
                \mat{Q}_{ji}                    & \text{ if } k=i      \\
                \mat{0} \in \R^{d_i \times d_k} & \text{ if } k \neq i
              \end{cases}
              \\
               & +
              \begin{cases}
                a_{(j,i)}(t) \mat{I}_{d_j}     & \text{ if } k=j      \\
                \mat{0}\in \R^{d_i \times d_k} & \text{ if } k \neq j
              \end{cases}
            \end{aligned}
            \label{eqn:M12_def}
          \end{align}
          %
          % The block \(\mat{M}_{12}(t)\) is given by \(\mat{M}_{12}(t) = \parens*{\mat{M}_{21}}^\transpose(t)\), due to the symmetric construction.\label{eqn:M21_def}
          %
          The block \(\mat{M}_{22}(t)\) comprises quadratic terms solely in \(\vec{e}\).
          As per the statement of the theorem, we express \(\mat{M}_{22}(t)\) as the sum of a block-diagonal matrix~\(\mat{\Lambda}\) containing the Poisson rates and a matrix \(\mat{R}(t)\) containing remaining terms.
          Thus, \(\mat{M}_{22}(t) = \mat{\Lambda} + \mat{R}(t)\), where
          \begin{subequations}\label{eqn:M22_def}%
            \begin{align}
               & \hphantom{\parens*{\mat{R}_2}_{((i, p), (j, q))}}
              \mathllap{\mat{\Lambda}} = \operatorname{diag}\parens*{ \ldots, \lambda_{(j,i)}\mat{I}_{d_j}, \ldots},
              \label{eqn:Lambda_def}
              \\
               & \parens*{\mat{R}}_{((j, i), (p, q))}(t) =
              \label{eqn:R2_def} % This label applies to the whole R_2 definition
              \\[1ex] % Line break added here
               & \qquad\qquad \begin{cases} % Indent the cases environment for readability
                                \parens*{ -2a_{(j,i)}(t) - \rho_{j} a_{(j,i)}(t)^2 } \mat{I}_{d_j} & \text{if } (j, i)=(p, q)             \\
                                -2\mat{Q}_{jp}                                                     & \text{if } i=q \text{ and } j \neq p \\
                                \mat{0}                                                            & \text{otherwise}
                              \end{cases}
              \notag
            \end{align}
          \end{subequations}
          The diagonal terms of \(\mat{R}(t)\) incorporate drift rates~\(a_{(j,i)}\) and weighting parameters~\(\rho_{j}\), while off-diagonal terms reflect cross-coupling from \(\mat{Q}\).

          The constant \(\gamma'\) aggregates all terms independent of \(\vec{s}\), arising solely from the upper bound of \(W_7\), where
          \begin{align}
            \gamma' = \sum_{i=1}^n \sum_{j \neq i} \tfrac{1}{\rho_{j}} \norm{\vec{y}_j^*}^2 
            % = \parens{n-1} \norm{\vec{y}^*}_{\mat{\mathrm{P}}^{-1}}^2.
            \label{eqn:gamma_prime_def}
          \end{align}
          We must show that for any given $\gamma > 0$ and $c_3 > 0$ (e.\,g., $c_3 = \beta$), constants $\rho_j > 0$ exist such that $\gamma' \le c_3 \gamma$.
          This condition is $\sum_{i=1}^{n}\sum_{j\ne i}\tfrac{1}{\rho_{j}}\norm{\vec{y}_{j}^{*}}^{2} \le c_3 \gamma$.
          As long as $\vec{y}^* \ne 0$ (if $y^*=0$, $\gamma'=0$ and \Cref{cor:mean-square-exponential-stability} applies), we can always find such $\rho_j$.
          For example, choosing $\rho_j = \rho$ for all $j$, the condition becomes $\rho \ge \tfrac{1}{c_3 \gamma} {(n-1)\norm{\vec{y}^{*}}^{2}}$ .
          Since the right-hand side is a finite positive constant, such $\rho_j$ always exist.
          We now proceed assuming such $\rho_j$ have been chosen.
          This fixes $\gamma'$ and also defines the matrix $\mat{R}(t)$ which depends on $\rho_j$.

    \item
          At any fixed instant~\(t\), \(\mat{M}(t)\) is a symmetric matrix.
          We apply \Cref{lemma:matrix-bounds} to this matrix for both statements of this theorem.
          \begin{subequations}
              \begin{enumerate}[
                label={\small(\alph*)},
                % wide
              ]
                    \item
                        We need to show \( \mat{M}(t) \succ \mat{0} \).
                        By \Cref{lemma:matrix-bounds}, this holds if all communication rates satisfy
                        \begin{equation}
                            \lambda_{(j,i)} > \lambda_{\text{s}}(t):= - \eigmin{ \mat{R}(t) - \mat{K}(t) },
                        \end{equation}
                        where \( \mat{K}(t) := \mat{M}_{21}^\transpose(t) \mat{M}_{11}^{-1} \mat{M}_{21}(t) \), and consequently, \( \mathcal{L}V(\vec{s}) \leq - \eigmin{\mat{M}(t)}\norm{\vec{s}}^2 + \gamma' \).
                        We set \({ c_3(t) = \eigmin{\mat{M}(t)} > 0}\).
                    \item
                        We need to show that \( \mat{M}(t) \succeq \beta \mat{I} \) for the given \( \beta \in \interval[open]{0}{2\eigmin{\mat{Q}}} \).
                        Let \( {\mu = \beta}\). 
                        Since \( {\mat{M}_{11} = 2 \mat{Q}} \), \(  \eigmin{ \mat{M}_{11} } = 2 \eigmin{\mat{Q}}  \) and the condition \( {\mu \in \interval[open]{0}{\eigmin{\mat{M}_{11}}}} \) from \Cref{lemma:matrix-bounds} is satisfied.
                        By \Cref{lemma:matrix-bounds}, \( \mat{M}(t) \succeq \beta \mat{I} \) holds, if all \( \lambda_{(j, i)} \) satisfy 
                        \begin{equation}
                            \lambda_{(j, i)} \geq \lambda_{\text{d}}(t) := \beta - \eigmin{\mat{R}(t) - \mat{K}_{\beta}(t)} ,
                        \end{equation}
                         where \( \mat{K}_{\beta}(t) := \mat{M}_{21}^\transpose(t) \parens{\mat{M}_{11} - \beta \mat{I} }^{-1} \mat{M}_{21}(t)\), and consequently \( \mathcal{L}V(\vec{s}) \leq - \beta \norm{\vec{s}}^2 + \gamma' \).
                        We set \( c_3 = \beta\).
              \end{enumerate}
          \end{subequations}

          To satisfy \Cref{lemma:gpues-ms-lyapunov-like-theorem} and guarantee stability for all time, we require a constant lower bound $c_3 \in (0,\infty)$ such that \({ \mathcal{L}V(\vec{s})\le -c_3\norm{\vec{s}}^{2}+\gamma'} \).
          This requires finding constant rates \( \lambda_{\text{s}} \) and  \( \lambda_{\text{d}} \) that satisfy the worst-case conditions \( \lambda_{(j,i)} > \lambda_{\text{s}} \ge \sup_{t}\{\lambda_{\text{s}}(t)\} \) (for (a)) and \( \lambda_{(j,i)} \ge \lambda_{\text{d}} \ge \sup_{t}\{\lambda_{\text{d}}(t)\} \) (for (b)).
          
          % The suprema are $\lambda_s = - \inf_t\{\mu_{min}(S(t))\}$ and $\lambda_d = \beta - \inf_t\{\mu_{min}(S_\beta(t))\}$ (where $S(t)=R(t)-K(t)$ and $S_\beta(t)=R(t)-K_\beta(t)$). 
          Finding the exact infimum of \( \eigmin{\cdot} \) over time is intractable.
          Instead, we construct constant matrices \( \mat{R}_c \) , \( \mat{K}_c \), and \( \mat{K}_{c,\beta} \) (hereafter $R$, $K$, $K_\beta$) that serve as uniform bounds, such that \( \mat{R}(t) \ge \mat{R} \), \( \mat{K}(t) \le \mat{K} \), and \( \mat{K}_\beta(t) \le \mat{K}_\beta \) for all $t$.
          
          This allows us to bound the infima, as \( \inf_t\{\eigmin{ R(t)-K(t) }\} \ge \eigmin{ \mat{R}- \mat{K}} \) and \( \inf_t\{\eigmin{R(t)-K_\beta(t)} \} \ge \eigmin{ \mat{R}- \mat{K}_\beta) } \).
          \begin{itemize}[wide]
              \item 
                Bounding \( \mat{R}(t) \): We require a constant matrix \( \mat{R} \) such that \( \mat{R}(t) \ge \mat{R} \).
                The diagonal entries of \( \mat{R}(t) \) are \(-2 a_{(j,i)}(t) -\rho_{j}a_{(j,i)}(t)^{2} \). 
                This is a downward-opening parabola in \( a_{(j,i)}(t) \).
                Given the assumption \( \abs{a_{(j,i)}(t) }\le a_{ji} \), its minimum on the interval \( \interval{a_{ji}}{a_{ji}} \) occurs at a boundary. 
                Thus, we define $\mat{R}$ as the constant matrix resulting from replacing \( a_{(j,i)}(t) \) with $a_{ji}$ in \eqref{eqn:R2_def}.
              \item
              Bounding \( \mat{K}(t) \) and \( \mat{K}_\beta(t) \):
              We require constant matrices \( \mat{K} \), \( \mat{K}_\beta \) such that \( \mat{K}(t) \le \mat{K} \) and \( \mat{K}_\beta(t) \le \mat{K}_\beta \) for all \( t \). 
              We find a valid upper bound \( \mat{K}_\beta \) by bounding its maximum eigenvalue \( \eigmax{ \mat{K}_\beta(t)} \le \norm{ \mat{K}_\beta(t)}_2 \).
              Using the submultiplicativity of matrix norms, we have  \( \norm{\mat{K}_\beta(t)}_2 \le \norm{ \mat{M}_{21}(t) }_2^2 \cdot \norm{\parens{ \mat{M}_{11}-\beta \mat{I}}^{-1} }_2 = \norm{ \mat{M}_{21}(t) }_2^2 \cdot \norm{ \parens{2\mat{Q}-\beta \mat{I}}^{-1}}_2 \).
              Next, we decompose \( \mat{M}_{21}(t) = \mat{M}_{21,c} + \mat{M}_{21,v}(t) \), where \( \mat{M}_{21,c} \) contains the constant \( \mat{Q} \)-terms and \( mat{M}_{21,v}(t) \) contains the time-dependent \( a_{(j,i)}(t) \) terms.
              We apply the triangle inequality \( \norm{ \mat{M}_{21}(t) }_2 \le \norm{ \mat{M}_{21,c} }_2 +  \norm{ \mat{M}_{21,v}(t) }_2 \).
              The norm of the constant part \( \norm{ \mat{M}_{21,c}}_2 \) is a fixed value.
              The norm of the (block-)diagonal variable part \( \mat{M}_{21,v}(t) \) is bounded by its largest entry in magnitude, using \( \abs{a_{(j,i)}(t)} \le a_{ji}\)  as \( \norm{ \mat{M}_{21,v}(t)}_2 = \max_{(j,i)}\{\abs{-a_{(j,i)}(t)}\} \le \max_{(j,i)}\{a_{ji}\} =: a_{\text{max}} \).
              Combining these steps, we obtain a constant upper bound \( {K}_{\text{bound},\beta} \) for the maximum eigenvalue of \( \mat{K}_\beta(t) \) as \( K_{\text{bound},\beta} := \parens{\norm{\mat{M}_{21,c}}_2 + a_{\text{max}}}^2 \cdot \norm{\parens{2 \mat{Q} - \beta \mat{I}}^{-1}}_2  \).
              We define \( \mat{K}_\beta = {K}_{\text{bound},\beta} \mat{I}\).
              This constant matrix satisfies \( \mat{K}_\beta(t) \le \eigmax{ \mat{K}_\beta(t)} \mat{I} \le \mat{K}_\beta \). 
              The identical procedure, setting \( \beta=0 \), yields the bound \( \mat{K} = K_{\text{bound}} \mat{I} \) for \( \mat{K}(t) \), where \( K_{\text{bound}} \) uses \( \norm{(2\mat{Q})^{-1}}_2 \).
          \end{itemize}
          With these constant bounds, we can define the sufficient constant rates
          \begin{subequations}
          \label{eqn:proof:thm:pes-ms:rates}
              \begin{align}
                \label{eqn:proof:thm:pes-ms:rates:s}
                \lambda_{\text{s}} &:=-\eigmin{ \mat{R}- \mat{K} }, \quad \text{and} \\
                \label{eqn:proof:thm:pes-ms:rates:d}
                \lambda_{\text{d}} &:= \beta - \eigmin{ \mat{R}- \mat{K}_\beta }.
              \end{align}
          \end{subequations}

    \item
        Finally, we verify the parameters of~\Cref{lemma:gpues-ms-lyapunov-like-theorem}.
        We chose \( c_1= \), \( c_2=1 \).
        In Step~2, we showed that for any given \( \gamma > 0 \) and \( c_3 > 0 \) (where \( c_3=\eigmin{\mat{M}_c } \) for (a) or  \( c_3=\beta \) for (b)), we can choose $\rho_j > 0$ such that the resulting  \( \gamma' = (n-1)\norm{ \vec{y}^{*} }_{P^{-1}}^{2} \) satisfies \( \gamma' \le c_3 \gamma \).
        The parameters from Lemma 1 are        \( \alpha = c_2/c_1 = 1 \) ,  \( \beta = c_3/c_2 = c_3 \) (which is \(\beta_{\text{Lemma}} \) for statement (b)), \( \gamma_{\text{Lemma}} = \gamma'/c_3 \le \gamma \).
        This confirms that the system is stable with the desired exponential rate (for (b)) and the desired ultimate bound $\gamma$.
    
          % Finally, \Cref{lemma:gpues-ms-lyapunov-like-theorem} relates he convergence rate and the region of ultimate convergence.
          % With \( c_3 = \eigmin{\mat{M}}\) and \(\gamma'\) as in \eqref{eqn:gamma_prime_def}, we conclude
          % \begin{subequations}
          %   \begin{align}
          %     \alpha & = 1
          %     \\
          %     \beta  & = \eigmin{\mat{M}}
          %     \\
          %     \gamma & = \frac{1}{\eigmin{\mat{M}}}\parens{n-1} \norm{\vec{y}^*}_{\mat{\mathrm{P}}^{-1}}^2.
          %   \end{align}
          % \end{subequations}
  \end{enumerate}
\end{proof}

}
 % \inversion{ecc2026}{\input{src/statements/thm.pes-ms.proof-sketch.tex}}
 \begin{proof}[of \Cref{cor:mean-square-exponential-stability}]
  Observe \(W_7 \equiv 0\) if \(a_{p \to i} = 0\) for all channels~\((p, i)\).
  Consequently, \(\gamma' = 0\) in the proof of \Cref{thm:practical-exponential-stability-ms}.
\end{proof}
\end{document}